\input amstex
\documentstyle{amsppt}
\NoBlackBoxes
\def\e{\varepsilon}
\def\a{\alpha}
\def\b{\beta}
\def\d{\delta}
\def\G{\Gamma}
\def\g{\gamma}
\def\k{\kappa}
\def\D{\Delta}
\def\l{\lambda}
\def\s{\sigma}
\def\x{\times}
\def \R{{I \!\! R}}
\def \N{{I \!\! N}}
\def \Z{{Z \!\!\! Z}}
\def \t{{T \!\!\! T}}
\def\o{\overline}
\def\f{\flushpar}
\def\u{\underline}
\def\v{\varphi}

\def\Om{\Omega}
\def\B{\Cal B}

\def\L{\Cal L}
\def\F{\Cal F}
\def\A{\Cal A}

\def\C{\Cal C}
\def\E{\Cal E}
\def\M{\Cal M}
\def\EE{\text{End}\,}

\def\({\biggl(}
\def\){\biggr)}
\def\lo{\longrightarrow}
\def\lra{\longrightarrow}
\def\dom{\Cal D}\document
\topmatter
\title
A salad of cocycles
\endtitle

\author{ Jon. Aaronson,\
 Mariusz Lema\'{n}czyk,\
  \&\
Dalibor Voln\'y .}\endauthor
\address{Aaronson:\ School of Mathematical Sciences, Tel Aviv University,
69978 Tel Aviv, Israel.}
\endaddress
\email{aaro\@math.tau.ac.il}
\endemail
\address{Lema\'{n}czyk: Institute of Mathematics,
 Nicholas  Copernicus University, \newline ul. Chopina 12/18, 87-100
Toru\'n, Poland.}
\endaddress
\email{mlem\@ipx1.mat.torun.edu.pl}
\endemail
\address{Voln\'y:\ Mathematical Institute, Charles University, Sokolovsk\'a 83,
186 00 Praha 8, Czech Republic}
\endaddress
\email{dvolny\@karlin.mff.cuni.cz}
\endemail
\thanks
Lema\'{n}czyk's research was partly supported by  KBN grant 2 P301
031 07.\smallskip Original uncut version of
\newline {\it A cut salad of cocycles.} Fund. Math. \bf{157},{\rm
(1998) 99-119.}
\endthanks
\abstract We study the centraliser of locally compact group
extensions of ergodic probability preserving transformations. New
methods establishing ergodicity of group extensions are
introduced, and new examples of squashable and non-coalescent
group extensions are constructed. Smooth versions of some of the
constructions are also given.
\endabstract
\endtopmatter
\heading
 \S0 Introduction
\endheading
Let $T$ be an  ergodic probability preserving
transformation of the probability space
$(X,\B,m)$.
\par Let $(G,\Cal T)$ be a locally compact, second countable,
topological  group ($\Cal T=\Cal T(G)$ denotes the family of open sets in
the topological space $G$), and let
$\v:X\to G$ be a measurable function.
\par The (left) {\it skew product} or $G$-{\it extension}
 $T_\v:X\x G\lo X\x G,$
is defined by
$$
T_\v(x,g)=(Tx,\v(x)g).
$$
The skew product preserves the measure $\mu=m\x m_G$ where $m_G$
is left Haar measure on $G$. There is an ergodic skew product
$T_\v:X\x G\lo X\x G$ iff the group $G$ is amenable (see
\cite{G-S}, references therein, and \cite{Zim}). In this paper, we
are mainly concerned with Abelian $G$. Recall that on any locally
compact, Abelian, second countable, topological  group $G$, there
is defined a {\it norm} $\|\cdot\|_G$ (satisfying $\|x\|=\|-x\|\ge
0$ with equality iff $x=0$,
 and $\|x+y\|\le\|x\|+\|y\|$) which generates the topology of $G$.

\

\subheading{The centraliser}
\par Recall that the {\it centraliser} of a non-singular transformation $R:X\to X$
is the collection of {\it commutors} of $R$, that is,
non-singular transformations of $X$ which commute with $R$.
The collection of invertible commutors (the {\it invertible
centraliser}) is denoted by $C(R)$.
\par We study those commutors $Q$ of $T_\v$, of form
$$Q(x,y)=(Sx,f(x)w(y))\tag{$*$}$$
where $w:G\lo G$  is a surjective, continuous group endomorphism, $S$
is a commutor of $T$,
and $f:X\lo G$ is measurable.
It is shown in proposition 1.1 of \cite{A-L-M-N} that
if $T$ is a Kronecker
transformation, and $T_\v$ is ergodic, then every commutor of $T_\v$
is of form ($*$).

\par Let $\EE(G)$ denote the collection of
surjective, continuous group endomorphisms of $G$ (a semigroup
under composition) and let
$$\E_\v=\{w\in\EE(G):\exists\text{ a commutor $Q$ of $T_\v$
of form ($*$) with } w=w_Q\},$$
a sub-semigroup of $\EE(G)$.

\

The study of $\E_\v$ yields  counterexamples:
\f if $\E_\v$ contains non-invertible endomorphisms, then $T_\v$ is
not {\it coalescent},
 i.e. its centraliser contains some non-invertible
transformation (see \cite{H-P}); and
\f if $\E_\v$ contains endomorphisms which do not preserve
$m_G$ (a possibility only for non-compact $G$),
then $T_\v$ is {\it squashable}, i.e.
 i.e. its centraliser contains some non-singular
transformation which is not measure preserving
 (see \cite{Aa1} and below).
Counterexamples like these (and others) will be discussed below.

\

\par In case $T$ is an odometer, for any
Abelian, locally compact, second countable $G$, the collection
$$\{\v:X\to G:T_\v\text{ ergodic },\ \E_\v=\{\text{Id}\}\}$$
is residual in the collection of measurable functions
$\v:X\to G$ considered in the topology
of convergence in measure (see below). Analogous results hold,
when $T$ is a rotation of the circle, for smooth $\v:\t\to\R$.

\

\subheading{Semigroup homomorphisms}
\par Let $\L_\v$ denote the collection of those commutors $S$ of $T$,
for which $\exists$ a commutor $Q$ of $T_\v$
of form ($*$) with $S=S_Q$.
\par When $G$ is Abelian and $T_\v$ is ergodic,
there is a surjective
semigroup homomorphism $\pi_\v:\L_\v\to\E_\v$ such that
if $S\in\L_\v,$ and $Q$ is a commutor of $T_\v$ of form ($*$) with
$S=S_Q$, then $w_Q=\pi_\v(S)$. This result (called
the {\it semigroup embedding lemma}) is proved at the end of this introduction.
\par It implies that $\E_\v$ is Abelian whenever the commutors of $T$
form an Abelian semigroup,
for instance when $T$ is a Kronecker transformation.
\par The restriction of $\pi_\v$ to
$L_\v(T)=\{S_Q:Q\in C(T_\v)$ of form $(*)\}$
is continuous with respect to the relevant Polish topologies by
the continuous embedding lemma established in \S1
(c.f. \cite{G-L-S} for the case where $G$ is compact).
\par The question arises as to when a homomorphism $\pi$ from a sub-semigroup $\Cal S$ of
commutors of $T$ into $\EE(G)$ occurs in this manner. That is, when does there exist a
measurable function $\v:X\to G$ such that $T_\v$ is ergodic, $\Cal S\subset
\L_\v$,
and $\pi=\pi_\v|_{ _{\Cal S}}$.
\par In \cite{L-L-T} it is shown that for
$T$ an invertible, ergodic probability preserving
transformation with some invertible commutor $S$ so that $\{S^mT^n:m,n\in\Z\}$ acts
freely, and $G=\t$, $\exists\ \v:X\to\t$ such that $S\in\L_\v,\ \E_\v\ni [x\mapsto 2x
\mod 1]$, and indeed, $\pi_\v(S)=[x\mapsto 2x \mod 1]$.
This includes the first example of a non-coalescent Anzai skew product
(i.e. $\t$-extension of a rotation of $\t$).

\

\subheading{The main results}
We generalise this to all Abelian, locally compact, second countable $G$:
\proclaim{Theorem 1} Suppose that $T$ is an ergodic
probability preserving transformation, $d\le\infty$,
and $S_1,\dots,S_d\in
C(T)\ (d\le\infty)$ are such that $(T,S_1,\dots,S_d)$ generate a free
$\Z^{d+1}$ action of probability preserving transformations of $X$.
\par If $w_1,\dots,w_d\in\EE(G)$ commute
 (i.e. $w_i\circ w_j=w_j\circ w_i\ \forall\ 1\le i,j\le d$), then there is
a measurable function $\v:X\to G$ such that
$$T_\v\text{ is ergodic,}$$
$$S_1,\dots,S_d\in\L_\v,\ \ w_1,\dots,w_d\in\E_\v;$$
and
$$\pi_\v(S_i)=w_i\ \ \ \ (1\le i\le d).$$
\endproclaim
Any Kronecker transformation of an uncountable compact group satisfies the
preconditions of theorem 1.
\proclaim{Theorem 2} Suppose that $T$ is an ergodic probability
preserving transformation,
and that $\{S_t:t\in\R\}\subset C(T)$ is such that $T$ and
$\{S_t:t\in\R\}$ generate a free
$\Z\x\R$ action of probability preserving transformations of $X$.
\par There is
a measurable function $\v:X\to \R$ such that
$$T_\v\text{ is ergodic;}$$
and there is a flow $\{Q_t:t\in\R\}\subset C(T_\v)$ of form
$$Q_t(x,y)=(S_tx,e^ty+\psi_t(x)).$$
In particular,
$$S_t\in\L_\v,\ \ w_t\in\E_\v\ \ \ \forall\ t\in\R$$
where $w_t(y)=e^ty$; and
$$\pi_\v(S_t)=w_t\ \ \ \ \forall\ t\in\R.$$
\endproclaim
\remark{Remarks}
\par 1) Theorem 1 can be extended (with analogous proof) to enable
"realisation" of a semigroup
homomorphism defined on a discrete subgroup of the centraliser which is
amenable, and which has F$\emptyset$lner sets which tile (see \cite{O-W}).
\par 2) In view of the reliance here on Rokhlin lemmas (see the proof of lemmas
4.1 and 4.2),
we ask if there is an ergodic probability preserving
transformation $(X,\B,m,T)$, and an ergodic $\v:X\to\t^2$ such that
$SL(2,\Z)\subset\E_\v$.
\par Note that $SL(2,\Z)\subset$ End$(\t^2)$, and that if $T$ is the $4$-shift
with symmetric product measure, then
$SL(2,\Z)\subset C(T)$.
\endremark
\subheading{Squashability and laws of large numbers}
\par Let $T=(X_T,\B_T,m_T)$ be a conservative, ergodic measure preserving
transformation of the $\s$-finite measure space $(X_T,\B_T,m_T)$.
Each commutor of $T$ is a measure multiplying transformation.
This is because for $Q$ a commutor of $T$, the measure $m_T\circ
Q^{-1}$ is
$m_T$-absolutely continuous, $T$-invariant, and hence
$${dm_T\circ Q^{-1}\over d\,m_T}\circ T={dm_T\circ Q^{-1}\over d\,m_T}$$
which is constant by ergodicity.
\par
The {\it dilation} of a measure multiplying transformation $Q$ is defined
by
$$D(Q)={dm\circ Q\over dm}\in(0,\infty].$$
Recall from \cite{Aa1} that the transformation
$T$ is called {\it squashable} if it has a commutor with non-unit dilation.
\par Let $\Cal C\subset\B_T,\ \Cal C=$ either $\B_T$ or
$\F_T=\{B\in\B_T:m_T(B)<\infty\}$. A {\it law of large numbers} for $T$
with respect to $\Cal C$
is a function $L:\{0,1\}^\N\to [0,\infty]$ such that
$$L(1_A,1_A\circ T,\dots)=m_T(A)\text{   a.e.  }$$
$\forall\ A\in\Cal C$. If $L$ is a law of large numbers for $T$
with respect to $\Cal C$, and $Q$ is a commutor of $T$ such that
$Q^{-1}\Cal C\subseteq\Cal C$, then $D(Q)=1$, as $\forall\ A\in\Cal C$,
$Q^{-1}A\in\Cal C$, and for a.e. $x\in X$,
$$m_T(A)=L(1_A(Qx),1_A(TQx),\dots)=L(1_{Q^{-1}A}(x),1_{Q^{-1}A}(Tx),\dots)
=m_T(Q^{-1}A).$$
Consequently,
\f {\sl if $T$ has  a law of large numbers with respect to
$\B_T$, then $T$ is non-squashable,}
\f and {\sl if $T$ has  a law of large numbers with respect to
$\F_T$, then no commutor of $T$ has non-unit, finite dilation.}
\par It was shown in \cite{Aa2,\ corollary 2.3, $\&$ theorem 3.4}
that {\sl if $G$ is a countable group
without arbitrarily large finite normal subgroups \f (e.g.
$G=\Z^k\x\Bbb Q^\ell$ or $G=\Z^\infty=\{(n_1,n_2,\dots)\in\Z^{\N}:n_k\to 0\}$),
\f then any ergodic
$G$-extension of a Kronecker transformation has  a law of large numbers with
respect to $\F$.}
\subheading{Example 1}
\par Let $T$ be a Kronecker transformation, then $\exists\ S\in C(T)$ so that
$\{S,\ T\}$ generate a free $\Z^2$ action. Let $G=\Z^\infty$,
and let $w=w_1\in\EE(G)$ be the shift $w((n_1,n_2,\dots))=(n_2,n_3,\dots)$.
By theorem 1, $\exists\ \v:X\to G$ such that $T_\v$ is ergodic, and $w\in\E_\v$.
\par It follows that:
\f {\sl $T_\v$ has a law of large numbers with respect to
$\F_{T_\v}$,}
\f but also a commutor $Q(x,y)=(Sx,w(y)+g(x))$ which has
infinite dilation since (as shown in the proof of proposition 1.1 of
\cite{A-L-M-N}) $D(Q)=D(w)=\infty$,
whence {\sl $T_\v$ has no law of large numbers with respect to $\B_{T_\v}$.}
\subheading{Complete squashability and Maharam transformations}
\par Evidently $D:C(T)\to\R_+$ is a multiplicative homomorphism.
Set
$\D_0(T)=D(C(T))$. The group $\D_0(T)$ was first considered in
\cite{H-I-K} (see also \cite{Aa2}).
If $T$ has a law of large numbers with respect to
$\F_T$, then $\D_0(T)=\{1\}$. In particular, (\cite{Aa2}, or
\cite{A-L-M-N}) if
$T$ is a $\Z$-extension of a Kronecker transformation, then
$\D_0(T)=\{1\}$. Our results on $\R$-extensions show that this result
fails dramatically for other transformations $T$, an $\R$-extension
of $T$ being a $\Z$-extension of a $\t$-extension of $T$. Moreover,
(see proposition 2.5) for Bernoulli $T$, any ergodic
$\R$-extension of $T$ is isomorphic to a $\Z$-extension of $T$.
\par  It is standard (see \S1 where we recall some well known facts about
Polish groups of measure multiplying transformations)
that $\D_0(T)$ is a Borel subgroup of $\R_+$.
The class
$$\{\D_0(T):T\text{ a conservative,
ergodic measure preserving transformation}\}$$
includes

\

\f {\sl $\R_+$, all countable subgroups of $\R_+$,} and
{\sl subgroups of $\R_+$
with any Hausdorff dimension} (see \cite{Aa2}).

\

\par Call a conservative, ergodic measure preserving transformation $T$ with
the property that $\D_0(T)=\R_+$ {\it completely
squashable}
Any ergodic {\it Maharam transformation} (defined below)
is completely squashable.
\f For a non-singular conservative, ergodic transformation $R$ of
$(\Om,\A,p)$,
the transformation $T:X=\Om\x\R\to X$ defined by
$$T(x,y)=(Rx,y-\log{dp\circ R\over d\,p})$$
preserves the measure $dm_T(x,y)=dp(x)e^ydy$, and
is called the {\it Maharam transformation} of $R$, as it was shown to be
conservative in \cite{Mah}.  If $Q_t(x,y)=(x,y+t)$, then $Q_t\in C(T)$
 and $D(Q_t)=e^t$.
\f Conservative, ergodic Maharam transformations were
constructed
in \cite{Kr}, and smooth Maharam transformations of $\t\x\R$ are
constructed in \cite{H-S}.
\par We show in \S5 that the transformations $T_\v$ constructed in
theorem 2 are isomorphic to Maharam
transformations (proposition 5.1), and
we obtain $\Z$-extensions of Bernoulli
transformations which are Maharam transformations (see the remarks after
proposition 5.1).
\par In \S-\S6 and 7, we present completely squashable
$\R$-extensions $T_\v$
which are not isomorphic to any Maharam transformation, for $T$ an
odometer, and $T$ a rotation of $\t$. For $T$ a rotation of $\t$,
our examples are as smooth as possible given the Diophantine properties
of the rotation number of $T$
(including  $C^\infty$, ergodic, completely squashable
$\R$-extensions for suitable rotation numbers).
\par It is not hard to construct   real analytic, ergodic,
 completely squashable
$\R$-extensions of a suitable irrational rotation using \S5 and
\cite{Kw-Le-Ru1} and \cite{Kw-Le-Ru2}.

\subheading{Conditions for ergodicity, and non-squashability of skew products}
\par Let $T$ be an  ergodic probability preserving
transformation of the probability space
$(X,\B,m)$, assume that $G$ is Abelian, and let $\v:X\to G$ be
measurable.
\par Recall from \cite{Sch}, that the {\it essential values} of $\v$ are defined by
$$E(\v)=\{a\in G:\forall\ A\in\B_+,\ a\in U\in\Cal T,\ \exists\ n\ge 1\ \ni\
m(A\cap T^{-n}A\cap [\v_n\in U])>0\},$$
which is a closed subgroup of $G$.
It is shown in \cite{Sch} that $T_\v$ is ergodic iff $E(\v)=G$.

\

Set
$$\widetilde D(\v)=\{a\in G:\exists\ q_n\ \ni\ T^{q_n}
\overset{\M(X)}\to\lra\text{Id},\ \&\ \v_{q_n}\to a\text{  a.e.}\}$$
where $\overset{\M(X)}\to\lra$ denotes convergence in the topology
of measure preserving transformations on $X$ (see \S1 below),
then (see \cite{A-L-M-N}) $E(\v)\supset$ Gp$\,(\widetilde D(\v))$ (the
group generated by $\widetilde D(\v)$).
\par If Gp$\,(\widetilde D(\v))$ is dense in $G$,
 then $T_\v$ is not only ergodic, but also
 non-squashable.

\

\par If $T$ is an odometer, then (see \cite{A-L-M-N}) for any
Abelian, locally compact, second countable $G$, there is a measurable
function $\v:X\to G$ such that Gp$\,(\widetilde D(\v))$ is dense in $G$.
Because of the density of coboundaries, such
functions are residual in  the collection of measurable functions
$\v:X\to G$ considered in the topology
of convergence in measure.
It is also shown in \cite{A-L-M-N} that there are rotations of $\t$
for which there is an ergodic, real analytic $\v:\t\to\R$ with
Gp$\tilde D(\v)$ dense in $\R$, whence
such functions are residual in any space (containing real analytic
functions) where the coboundaries are dense.

\

\par It is shown in \cite{L-V}, that  for a rotation $T$ of $\t$,
there exist
dense $G_\delta$ sets in the spaces of absolutely continuous, Lipschitz,
$k$ times
continuously or infinitely differentiable functions $f$ with zero mean on
$\Bbb T$ for which the distributions of $f_{n_k}$ converge to a continuous
distribution along a rigid sequence $n_k\to \infty$ whenever theses spaces
contain  non-trivial cocycles (i.e. not $T$-cohomologous to a constants).
\par For irrational rotations
with bounded partial quotients, nontrivial cocycles
exist in the space of
absolutely continuous functions (while every zero mean Lipschitz function is
a coboundary). For rotations with unbounded partial quotients,
nontrivial cocycles exist in the space $\Cal C^p$ of $p$ times
continuously differentiable functions if and only if
\f $\limsup_{n\to\infty} q_{n+1}/q_n^p=\infty$
where $\{q_n:n\in\N\}$ are the principal denominators of the rotation
(cf. \cite{Ba-Me}).
\f For irrational rotations satisfying this condition $\forall\
p\in\N$, there are non-trivial infinitely differentiable cocycles.

\par Every such cocycle (the distributions of which converge to a continuous
one along a rigid sequence) is ergodic (see \cite{L-V}).
But also $\E_f\subset\{\pm 1\}$, since, if
$f\circ S= cf+g\circ T-g$, where $S$ is another rotation of $\t$;
\f on the one hand, the distributions $f_{n_k}\circ S$
converge to the same limit as the distributions of $f_{n_k}$,
\f while on the other hand, $(g-g\circ T)_{n_k}\to 0$
in measure, hence the limit distribution is invariant under multiplication
by $c$ which implies $c=\pm 1$ and $T_f$ non-squashable.

\

\subheading{New conditions for ergodicity allowing squashability}
The conditions for ergodicity of skew products discussed in
\cite{A-L-M-N} and \cite{L-V} are unsuitable for our constructions of
squashable skew products as they eliminate squashability.

We need  new conditions for the
ergodicity of a measurable function $\v:X\to G$ which are flexible
enough to allow squashability.

Such conditions, called {\it essential value conditions} are introduced in \S3.

Cocycles are constructed as infinite sums of coboundaries. Each coboundary
"contributes" a particular essential value condition, which the subsequent
coboundaries are "too small" to destroy. The
essential value conditions remaining for the infinite sum gives
its ergodicity.

The simplest version of our  essential value conditions
 is the {\it rigid } one (see proposition 6.2)
used in the constructions of \S6  and \S7 which could form an introduction
to the proofs of theorems 1 and 2 in \S4.

\

To conclude this introduction, we prove the
\proclaim{Semigroup embedding lemma}Suppose that $G$ is Abelian, and
that $\v:X\to G$ is such that $T_\v$ is ergodic. There is a
surjective semigroup homomorphism
$$\pi_\v:\L_\v\to\Cal E_\v$$
such that
if $Q(x,y)=(Sx,f(x)w(y))$ defines a commutor of $T_\v$,
then $w=\pi_\v(S)$.
\endproclaim
\demo{Proof}We must show that if $S\in\L_\v,\ w_1,\ w_2\in\E(G)$,
$f_i:X\to G,\ (i=1,2)$ are measurable, and
$Q_i(x,y)=(Sx,f_i(x)w_i(y))$ are such that
$Q_i\circ T_\v=T_\v\circ Q_i,\ (i=1,2)$, then
$w_1=w_2$.
\par To this end, let $U=w_1-w_2$, then $T_{U\circ\v}$ is an ergodic
transformation of $X\x U(G)$ (being a factor of $T_\v$ via
$(x,y)\mapsto (x,U(y))$).
The condition $Q_i\circ T_\v=T_\v\circ Q_i$ means that
$$\v\circ S=w_i\circ\v+f_i\circ T- f_i\ \ \ (i=1,2),$$
whence
$$U\circ\v=g\circ T-g$$
where $g=f_1-f_2$.
Define $\tilde g:X\to G/U(G)$ by $\tilde g(x)=g(x)+U(G)$.
It follows that $\tilde g\circ T=\tilde g$, whence by ergodicity of $T$,
$\exists\ \g\in G$ such that $\tilde g=\g+U(G)$ a.e. Therefore
$h:=g-\g:X\to U(G)$ is measurable and satisfies
$$U\circ\v=h\circ T-h.$$
The ergodicity of $T_{U\circ\v}$ on $X\x U(G)$ now implies $U(G)=\{0\}$,
i.e. $U\equiv 0$, or $w_1=w_2$.
\par We've shown that $\forall\ S\in\L_\v,\ \exists !\ w=:\pi_\v(S)\in\E_\v$
such that
$\exists\ f_S:X\to G$ measurable so that
$Q(x,y)=(Sx,f_S(x)\pi_\v(S)(y))$ defines a commutor of $T_\v$.
The rest of the lemma follows easily from this.
\hfill\qed\enddemo

\newpage
\heading
\S1 polish groups of
measure multiplying transformations
\endheading
\par Given a $\s$-finite measure space $(Y,\C,\nu)$ let
$\Cal M(Y,\C,\nu)$ denote the group of invertible measure multiplying
transformations of $(Y,\C,\nu)$, i.e. non-singular transformations $Q:Y\to Y$
such that $\nu\circ Q^{-1}=c\nu$ for some constant $c\in\R_+$.
This is a Polish group when equipped with the weak topology inherited from
that of the invertible, bounded linear operators on $L^2(Y,\C,\nu)$.
A metric for this topology is defined by
$$\rho(Q,R)=\sum_{n=1}^\infty \(\|f_n\circ Q-f_n\circ R\|_2+\|f_n\circ Q^{-1}-f_n\circ R^{-1}\|_2\)$$
where $\{f_n:n\in\N\}$ is a C.O.N.S. in $L^2(m_G)$.
The dilation function $D:\M\to\R_+$ defined (as above)
by $Q\mapsto{d\nu\circ Q^{-1}\over d\nu}:=D(Q)$ is a continuous homomorphism.
\par Let $T$ be a conservative, ergodic measure preserving transformation of
the standard $\s$-finite measure space $(X,\B,m)$, then $C(T)$ is
a closed subgroup of \f$\Cal M(X,\B,m)$ and hence Polish. The
multiplicative homomorphism $D:C(T)\to\R$ is continuous, and so
Ker$\,D$ is a closed, normal subgroup of $C(T)$. The natural
topology on $C(T)/$Ker$\,D$ is Polish, and $D:C(T)/$Ker$\,D\to\R$
is continuous and injective. By Souslin's theorem (see \cite{Kur,\
\S36,\ IV}),
$$\D_0(T)
=D(C(T)/\text{Ker}\,D)$$ is a Borel set in $\R$.
\par Let $G$ be a locally compact, second countable
topological  group, then $G$ is $\s$-compact, and Polish. Let $m_G$ be left
Haar measure on $G$.
\par The action of $G$ on $(G,\B(G),m_G)$ by left
translation is ergodic. To see this, suppose $A\in\B(G)_+$ and
$m_G(gA\D A)=0\ \forall\ g\in G$. The measure $m'$
defined by $dm'=1_Adm_G$ is a left Haar measure on $G$, and by
unicity of such, $m'=m_G$ whence $A=G$ mod $m_G$. The ergodicity of the
action of $G$ by right translation is obtained in a similar manner, as
right Haar measure is equivalent to $m_G$.
\par The maps $L,\ R:G\to\M(G)$ given by $L_gf(x):=f(gx),\ R_gf(x):=f(xg)$
are continuous, and their ranges are closed in $\M(G)$. This follows
from the ergodicity of the actions of $G$ by translation; indeed, if
$R=\lim_{n\to\infty}R_{g_n},$ set $f(x)=x^{-1}R(x)$,
which is $L_g$-invariant $\forall\ g\in G$, hence constant, and
$R=R_h$ for some $h\in G$. Let the ranges of these maps in $\M(G)$ be
$\tilde G_L$ and $\tilde G_R$, considered with their inherited (Polish)
topologies. By Souslin's theorem (\cite{Kur,\ \S36,\ IV}), the inverse maps $L^{-1}:\tilde G_L\to G$
and $R^{-1}:\tilde G_R\to G$ are both measurable, and (being group isomorphisms)
are continuous by Banach's theorem (\cite{Ban,\ p.20}). In particular, the metric $d$ on $G$, defined by
$d(x,y)=\rho(L_x,L_y)$, generates the topology of $G$.
\par In case $G$ is Abelian,
$\|x\|_G:=d(L_x,\text{Id})$ defines a topology generating norm on $G$.
\par Let $\text{Aut}(G)$ denote the group of
continuous group automorphisms of $G$. For example, $\tau_g(x):=g^{-1}xg$
is a continuous group automorphisms of $G$ (called an {\it inner} automorphism).
Also, consider Aff$(G)$, the group of invertible, {\it affine } transformations of $G$
of form $L_g\circ w\equiv R_g\circ w'$ where $g\in G$ and $w,w'\in\text{Aut}(G)$
(and $w'=\tau_g\circ w$).
\par Both Aff$(G)$ and $\text{Aut}(G)$ are closed subgroups of $\M(G)$.
To show that Aff$(G)$ is closed in $\M(G)$, we note first that
if $Q\in\M(G)$, then $Q\in$Aff$(G)$ iff $\exists\ w\in\text{Aut}(G)$ such that
$$Q\circ L_g=L_{w(g)}\circ Q\ \ \forall\ g\in G.$$
Indeed, supposing this condition, the function $x\mapsto w(x)^{-1}Q(x)$
is $L_g$-invariant $\forall\ g\in G$, hence constant, and
$Q=R_h\circ w\in$ for some $h\in G$.
\par Also, the topology on Aff$(G)$
inherited from $\M(G)$ coincides with the compact-open topology. This follows because
Aff$(G)$, equipped with the compact-open topology, is a Polish space continuously embedded
(by the identity) onto Aff$(G)$ equipped with the topology inherited from $\M(G)$;
this identity necessarily being a homeomorphism (as shown above
by the theorems of Banach and Souslin).
\par It is not hard to show that the Polish topology on  Aff$(G)$
also coincides with the topology of pointwise convergence.
\par It follows from the above that Aut$(G)$ is closed in Aff$(G)$.
We'll assume throughout that Inn$(G):=\{\tau_g:g\in G\}$
is also closed in Aff$(G)$. As pointed out to us by Danilenko and Glasner,
this is the case e.g. when $G$ is
Abelian, compact or a connected Lie group; but not in general.
Let $T$ be an ergodic probability preserving transformation of the
standard non-atomic probability space
$(X,\B,m)$, and let
$G$ a locally compact, second countable,
topological  group with left Haar measure $m_G$, and
consider the $\s$-finite measure space $(X\x G,\B\x\B(G),\mu)$ where $\mu=m\x m_G$.

The main point of this section is to establish the continuous embedding
lemma (see below) which is a topological version of the semigroup
embedding lemma.
Let $\Cal M=\Cal M(X\x G,\B\x\B(G),\mu)$, and
let $\widetilde{\M}$ denote those
$Q\in\Cal M$ of form
$$Q(x,y)=(Sx,h(x)w(y))\tag{$*$}$$
where $w\in\text{Aut}(G)$  is a continuous group automorphism, $S\in \M(X)$
and $h:X\lo G$ is measurable.

\

Write
$$S_f(x,y)=(Sx,f(x)y)$$
for $S\in\M(X)$ and $f:X\to G$ measurable.
Also, for $w\in\text{Aut}(G)$, write
$$W_w(x,y)=(x,w(y)).$$
\par If $Q\in\widetilde{\M}$ is as in ($*$),
then
$$Q=S_{f}\circ W_{w}.$$
It is clear that the representation $S_f(x,y)=(Sx,f(x)y)$ is unique in the sense
that $S_f=S'_{f'}$ implies $S=S'$ and $f=f'$. The unicity of the representation of $Q\in\widetilde{\M}$
by ($*$) follows from this, and $\s_{w(g)}=Q\circ\s_g\circ Q^{-1}$ where
$\s_g(x,y)=(x,yg)$.

\

Now let $T:X\to X$ be an ergodic probability preserving transformation, and let $\v:X\to G$ be
measurable.
\par We study $\widetilde C(T_\v):=C(T_\v)\cap\widetilde{\M}$
- a closed subgroup of $\M$.
If $Q\in\widetilde{\M}$, then
$Q\in C(T_\v)$
iff $S_Q\in C(T)$, and
$$\v\circ S_Q=h_Q\circ T\cdot w_Q\circ\v\cdot h_Q^{-1}.$$

\

Set $$\A_\v(T)=\{w_Q:Q\in\widetilde C(T_\v)\},\ \&\ L_\v(T)=\{S_Q:Q\in\widetilde C(T_\v)\}=\L_\v\cap C(T).$$
\proclaim{Continuous embedding lemma}
Suppose that $T$ is an ergodic probability preserving transformation of the standard probability space
$(X,\B,m)$, that $G$ is a locally compact, second countable topological group, and that
$\v:X\to G$ is measurable such that $T_\v$ is ergodic.
\par There is a Polish topology on $L_\v$, stronger than the topology inherited from
$\M(X)$, and a continuous homomorphism
$$\pi_\v:L_\v\to \text{Aut}(G)/\text{Inn}(G)$$
such that if $Q\in C(T_\v)$ is of form $(*)$, then
$$w_Q\text{Inn}(G)=\pi_\v(S_Q).$$
In case $G$ is Abelian, then $\pi_\v$ is the restriction to $L_\v$ of the homomorphism in the
semigroup embedding lemma, and
there is a Polish topology on $\A_\v(T)$ stronger
than that
inherited from Aut$(G)$ such that
$$\pi_\v:L_\v(T)\to \A_\v(T)$$
is continuous.
\endproclaim
The proof of the continuous embedding lemma uses four lemmas, two of
which concern the structure of $\widetilde{\M}$.
\par Let $\o G:=\{\s_g:g\in G\}$, then $\o G=\{$Id$\}\x\tilde G_R$ is
a closed subgroup
of $\M$, and the embedding $g\mapsto\s_g$ is a homeomorphism
( $G\leftrightarrow \o G$). Also, $\o G\subset\widetilde{\M}$ because
$\s_g(x,y)=(x,g(g^{-1}yg))$. We'll need
\proclaim{Lemma 1.1}
If $Z$ is a separable metric space, and
$f:X\x G\to Z$ is measurable and $f\circ\s_g=f$ a.e. $\forall\ g\in G$,
then $\exists\ g:X\to Z$ such that $f(x,y)=g(x)$ a.e.
\endproclaim
\demo{Proof}
Choose $h:Z\to [0,1]$ injective, and (Borel) measurable. For $A\in\B$, the
function $f_A:G\to\R$ defined by
$$f_A(y)=\int_Ah(f(x,y))dm(x)$$
is $R_g$-invariant $\forall\ g\in G$, and hence, by
ergodicity of $\tilde G_R$ on $G$, $\exists\ c(A)\in\R$ such that
$f_A=c(A)$ a.e. Since $c:\B\to\R$ is a $m$-absolutely continuous
signed measure, $\exists\ k:X\to\R$ such that $h(f(x,y))=k(x)$ a.e.,
and the required function is $g(x)=h^{-1}(k(x))$.
\hfill\qed\enddemo
Suppose that
$Q\in\widetilde{\M}$, then, as mentioned above
$$\forall\ g\in G,\ \exists\ g'\in G\ \ni\ Q\circ\s_g=\s_{g'}\circ Q.$$
We obtain the converse statement as an immediate consequence of lemma 1.1.
\par Supposing that $Q\in\M$ and that
$$\forall\ g\in G,\ \exists\ g'\in G\ \ni\ Q\circ\s_g\circ Q^{-1}=\s_{g'},$$
we obtain a continuous group endomorphism $w:G\to G$ such that
$$Q\circ\s_g\circ Q^{-1}=\s_{w(g)}.$$
Writing
$$Q(x,y)=(S(x,y),F(x,y)),$$
we have that
$$S(x,yg)=S(x,y),\ \ \ F(x,yg)=F(x,y)w(g).$$
The functions $(x,y)\mapsto S(x,y)$, and $(x,y)\mapsto F(x,y)w(y)^{-1}$
are $\s_g$-invariant $\forall\ g\in G$, and hence, by
lemma 1.1, for a.e. $(x,y)\in X\x G$
$$S(x,y)=S(x)\text{  and }F(x,y)=f(x)w(y).$$
The assumption that $Q\in\M$ now gives that
$Q\in\widetilde{\M}$.

\

We now discuss the topology of $\widetilde{\M}$.
It can easily be shown that the topology inherited by
$$\{\text{Id}_f:f:X\to G\text{  measurable}\}$$
from $\M$ is the topology of convergence in measure.
\proclaim{Lemma  1.2}
$\widetilde{\M}$ is closed in $\M$, and the projections
$$Q\mapsto S_Q,\ Q\mapsto f_Q\text{  and }Q\mapsto w_Q$$
are continuous.
\endproclaim
\demo{Proof}Firstly, suppose that $Q_n\in\widetilde{\M}$, and
$Q_n\to Q$ in $\M$.
We have that
$\s_{w_n(g)}=Q_n\circ\s_g\circ Q_n^{-1}$ converges, necessarily to
$Q\circ\s_g\circ Q^{-1}=\s_{w_(g)}$ since $\o G$ is closed in $\M$, and
by the above application of lemma 1.1, $Q\in\widetilde{\M}$. This also proves that
$$w_{Q_n}^{\pm 1}\to w_Q^{\pm 1}\text{ pointwise,}$$
and hence in $\text{Aut}(G)$.

\

To see that $Q\mapsto S_Q$ is continuous, let $Q,R\in\widetilde\M$.
Let $A\in\B(X)$, then $\exists\ C\subset A,\ C\in\B(X)$ such that
$m(C)>{m(S_Q^{-1}A\D S_R^{-1}A)\over 2}$, and
\f either (a) $S_Q^{-1}C\cap S_R^{-1}A=\emptyset$, or (b)
$S_Q^{-1}A\cap S_R^{-1}C=\emptyset$.
\par Now let $F\in\B(G),\ m_G(F)<\infty$.
\f In case (a), $Q^{-1}(C\x F)\cap R^{-1}(A\x F)=\emptyset$, and
$$\mu(Q^{-1}(A\x F)\D R^{-1}(A\x F))\ge \mu(Q^{-1}(C\x F))=D(Q)m(C)m_G(F),$$
and similarly,in case (b),
$Q^{-1}(A\x F)\cap R^{-1}(C\x F)=\emptyset$, and
$$\mu(Q^{-1}(A\x F)\D R^{-1}(A\x F))\ge \mu(R^{-1}(C\x F))=D(R)m(C)m_G(F).$$
This shows that
$$m(S_Q^{-1}A\D S_R^{-1}A)\le
{2\mu(Q^{-1}(A\x F)\D R^{-1}(A\x F))\over \min\{D(Q),D(R)\}m_G(F)}.$$

The continuity of $Q\mapsto D(Q)=D(w_Q)$ now shows the continuity of
$Q\mapsto S_Q$.
\par Finally, the continuity $Q\mapsto f_Q$
follows from
$$(\text{Id})_{f_Q}=(S_Q)_0^{-1}\circ Q\circ W_{w_Q}^{-1}.$$
\hfill\qed
\enddemo
\proclaim{Lemma 1.3}
\par There is a topology on $L_\v(T)$ with respect to which it
is a Polish group, continuously embedded in $C(T)$, and
$$L_\v(T)\cong\widetilde C(T_\v)/\o G.$$
\endproclaim
\demo{Proof} Write $p(Q)=S_Q$, then $p:\widetilde{\M}\to\M(X)$ is continuous, and
$p(\widetilde C(T_\v))=L_\v(T).$ Clearly, $p\circ \s_g=p$ for all $g\in G$,
whence $p:\widetilde C(T_\v)/\o G\to L_\v(T)$ is well-defined, onto, and
continuous when $\widetilde C(T_\v)/\o G$ is equipped with the quotient
(Polish) topology.
We claim that $p|_{\widetilde C(T_\v)/\o G}$ is actually injective.
\par To see this, suppose that
$Q\in\widetilde C(T_\v)$ and $Q(x,y)=(x,h(x)w(y))$, then
$$h\circ T\cdot w\circ\v=\v\cdot h.$$
Now set $F(x,y)=y^{-1}h(x)w(y)$, then
$$\align F\circ T_\v(x,y) &=y^{-1}\v(x)^{-1}[h(Tx)w(\v(x))]w(y)\\ &=
y^{-1}\v(x)^{-1}[\v(x)h(x)]w(y)\\ &=
F(x,y)\endalign$$
and $F$ is constant by ergodicity of $T_\v$, whence
$Q\in\o G.$
\par
The group isomorphism $p:\widetilde C(T_\v)/\o G\to L_\v(T)$ can be
used to transport the Polish structure to $L_\v(T)$ which is
a Polish group, continuously embedded in $C(T)$.
\hfill\qed\enddemo
\remark{Remark}
By Souslin's theorem (see \cite{Kur,\ \S36,\ IV}:
$L_\v(T)$ is a Borel subset of $C(T)$, and
$$\B(L_\v(T))=\B(C(T))\cap L_\v(T).$$
\endremark

Let
$\widetilde C_0(T_\v)=\{Q\in\widetilde C(T_\v):w_Q=\text{Id}\}$,
a normal, closed subgroup of $\widetilde C(T_\v)$.
\proclaim{Lemma 1.4}
\par There is a topology on $\A_\v(T):=\E_\v\cap$Aut$(G)$
with respect to which it
is a Polish group, continuously embedded in
$\text{Aut}(G)$,
and, as Polish groups,
$$\A_\v(T)\cong\widetilde C(T_\v)/\widetilde C_0(T_\v).$$
\endproclaim
\demo{Proof} Write $q(Q)=w_Q$ for $Q\in \widetilde C(T_\v)$. By
lemma 1.2, $q:\widetilde C(T_\v)\to\text{Aut}(G)$ is continuous, and
$q(\widetilde C(T_\v))=\A_\v(T)$.
Clearly,
$$\text{Ker }q=\widetilde C_0(T_\v),$$
whence $q:\widetilde C(T_\v)/\widetilde C_0(T_\v)\to\A_\v(T)$
is well-defined and bijective.
\par As before, the group isomorphism can be
used to transport the quotient Polish topology on
$\widetilde C(T_\v)/\widetilde C_0(T_\v)$ to $\A_\v(T)$ which thus becomes
a Polish group, continuously embedded in $\text{Aut}(G)$.
\hfill\qed\enddemo
\demo{Proof of the continuous embedding lemma}
Let
$$\widetilde C_I(T_\v)=\{Q\in\widetilde C(T_\v):w_Q\in\text{Inn}(G)\}.$$
Note that $\widetilde C_I(T_\v)$ is a closed normal subgroup of $\M$, and is
generated by $\widetilde C_0(T_\v)$, and $\o G$.

It now follows that
\par $\widetilde C(T_\v)/\widetilde C_I(T_\v)$ is a Polish group, continuously
embedded in $\text{Aut}(G)/\text{Inn}(G)$ by
$$Q\widetilde C_I(T_\v)\mapsto w_Q\text{Inn}(G).$$
A natural map $L_\v(T)\to\text{Aut}(G)/\text{Inn}(G)$ is now generated by
$$L_\v(T)\cong \widetilde C(T_\v)/\o G\to\widetilde C(T_\v)/\widetilde
C_I(T_\v)\to\text{Aut}(G)/\text{Inn}(G).$$
It follows from the above that this map is continuous.

\

In case $G$ is Abelian, $\text{Inn}(G)=\{$Id$\}$, and the above
becomes a statement of the continuity of:
$$L_\v(T)\cong \widetilde C(T_\v)/\o G\to\widetilde C(T_\v)/\widetilde
C_0(T_\v)\cong\A_\v(T).$$
\hfill\qed\enddemo

\newpage
\heading
 \S2 Eigenvalues of skew products
\endheading
Let $T$ be an  ergodic probability preserving
transformation of the probability space
$(X,\B,m)$,
let $G$ be an locally compact, second countable,
topological  group, and let
$\v:X\to G$ be a cocycle with
$T_\v$ ergodic on
$X\times G.$

\

Recall that if $R$ is non-singular, conservative, ergodic and
$f:X\to\Bbb C$ is measurable such that
$$f\circ R=\l f,$$
where $\l\in\Bbb C$, then $|f|$ is constant (w.l.o.g. $=1$), $|\l|=1$.
\par We consider the situation where the measurable function $\v:X\to G$
is {\it aperiodic} in the sense that
all eigenvalues for the skew product $T_\v$
are eigenvalues for $T$ (that is, if $f:X\x G\to\t$ is measurable and
$f\circ T_\v=\l f$ where $\l\in\t$, then $\exists\ g:X\to\t$
measurable such that
$f(x,y)=g(x)$ a.e.).
The main result of this section is
\proclaim{Proposition 2.1}
\f If $G=\R$ or $\t$, $T_\v$ is ergodic, and $\E_\v\ne\{$Id$\}$,
then $\v$ is aperiodic.\endproclaim
\proclaim{Lemma 2.2}
\f Suppose that $T_\v$ is ergodic and $f:X\x G\to\t$ is measurable
such that $f\circ T_\v=\l_0 f$ where $\l_0\in\t$, then
$$f(x,y)=f_0(x)\g(y)\text{ where }\g\in\widehat G\text{ and }f_0:X\to\t\text{ is measurable.}$$
\endproclaim
\remark{Remark} Note that we do not assume that $G$ is Abelian here.
\endremark
\demo{Proof} For $Q\in C(T_\v)$, we have that
$$(f\circ Q)\circ T_\v=f\circ T_\v\circ Q=\l_0 f\circ Q,$$
whence, by ergodicity of $T_\v$,
$\exists\ \l(Q)\in\t$ such that $f\circ Q=\l(Q) f$ (note that $\l(T_\v)=\l_0$.  The mapping
$\l(Q):C(T_\v)\to\t$ is a continuous homomorphism.
\par Since $\o G\subset C(T_\v)$, we obtain $\g\in\widehat G$ by taking
$\g(g):=\l(\s_g)$. Thus
$$f\circ \s_g=\g(g)f\ \ \forall\ g\in G.$$
Set $F(x,y)=\g(y)^{-1}f(x,y)$, then $F\circ \s_g=F\ \forall\ g\in G$, whence
by lemma 1.1, for a.e. fixed $x\in X,\ F(x,\cdot)$ is constant.
This proves the lemma.
\hfill\qed\enddemo
\proclaim{Lemma 2.3}
\f Suppose that $T_\v$ is ergodic and $f=f_0\otimes\mu$ where
$f_0:X\to\t$ is measurable, $\g\in\widehat G$, and
$f\circ T_\v=\l_0 f$ for some $\l_0\in\t$, then
$$\g\circ w=\g\ \ \forall\ w\in\E_\v.$$
\endproclaim
\demo{Proof}Let $\l:C(T_\v)\to\t$ be such that
$$f\circ Q=\l(Q) f\ \ \forall\ Q\in C(T_\v).$$
Suppose that $w\in\E_\v$, and $Q\in\widetilde C(T_\v)$ with
$Q(x,y)=(Sx,h(x)w(y))$, then
$$\align \l(Q)f_0\otimes\g(x,y) &=\l(Q)f(x,y)\\
& = f\circ Q(x,y)\\ &= f_0(Sx)\g(h(x))\g(w(y))\\ & =
[(f_0\circ S)\cdot(\g\circ h)]\otimes\g\circ w (x,y),\endalign$$
and since the character $\g\in\widehat G$ appearing in the
eigenfunction $f_0\otimes\g$ is unique,
$$\g\circ w=\g.$$
\hfill\qed\enddemo
\demo{Proof of Proposition 2.1}
This now follows from lemma 2.3, as
if $G=\t,\ \R$, and $\g\in\widehat G,\ w\in$End$(G)$, then
$\g\circ w=\g$ iff either
$\g\equiv 1$ or $w=$Id.
\hfill\qed\enddemo
\proclaim{Corollary 2.4 \ \ \ (cf \cite{Rob})}
\f If $T_\v$ is ergodic, then the only eigenvalues of
$T_\v$ are the eigenvalues of $T_{\v[G,G]}:X\x G/[G,G]
\to X\x G/[G,G].$
\endproclaim
\demo{Proof}
By lemma 2.2, an eigenfunction of $T_\v$ must be of
form
$f_0\otimes\mu$ where $\mu\in\widehat G$. But
$$\widehat G=\widehat{G/[G,G]},$$
so that any eigenfunction of $T_\v$ is actually an eigenfunction of
$T_{\v[G,G]}$.
\hfill\qed\enddemo
\proclaim{Proposition 2.5}
\f If $T$
is Bernoulli, then any ergodic $\R$-extension of $T$ is (isomorphic to)
a $\Z$-extension of $T$.
\endproclaim
\demo{Proof} Let $\v:X\to\R$ be an measurable such that $T_\v$ is ergodic.
For $c>0$ let $\v^{(c)}:X\to \t\cong [0,c)$ be defined by $\v^{(c)}=\v$ mod $c$.
For each  $c>0$ there is a measurable function $\psi^{(c)}:X\x\t\to\Z$ such that
$$T_\v\cong (T_{\v^{(c)}})_{\psi^{(c)}}.$$
It is known that for some $c>0$, $\v^{(c)}$ is not cohomologous to
a constant, else (see \cite{M-S} and \cite{H-O-O})
$\v$ would be cohomologous to
a constant and not ergodic. For this $c>0$, $T_{\v^{(c)}}$ is weakly mixing, whence by theorem 1 of \cite{Rud} $T_{\v^{(c)}}$ is Bernoulli, and since $h(T_{\v^{(c)}})=h(T),$ we have by \cite{Or} that $T_{\v^{(c)}}\cong T$. The conclusion is that $(T_{\v^{(c)}})_{\psi^{(c)}}$ is a $\Z$-extension of $T$.  \hfill\qed\enddemo \heading \S3 essential value conditions\endheading
Let $T$ be an invertible, ergodic probability preserving transformation of the standard probability space
$(X,\B,m)$, let $G$ be a locally compact, second countable Abelian group, and let $\v:X\to G$
be measurable. We develop here a countable condition for ergodicity of
$T_\v$.
The EVC's to be defined are best understood in terms
of orbit cocycles, and the {\it groupoid} of $T$ (see \cite{Fe-Mo}).
\par A {\it partial probability preserving transformation} of $X$ is a pair $(R,A)$ where $A\in\B$ and $R:A\to RA$ is invertible and
$m|_{RA}\circ R^{-1}=m|_A$. The set $A$ is called the {\it domain } of $(R,A)$.
We'll sometimes abuse this notation by writing $R=(R,A)$ and $A=\dom(R)$.
Similarly, the {\it image} of $(R,A)$ is the set $\Im(R)=RA$.

\par The {\it equivalence relation} generated by $T$ is
$$\Cal R=\{(x,T^nx):x\in X,\ n\in\Z\}.$$
\par For $A\in\B(X)$ and $\phi:A\to\Z$, define $T^\phi:A\to X$ by
$T^\phi(x):=T^{\phi(x)}x$.
\par The {\it groupoid} of $T$ is
$$[T]=\{T^\phi:\ T^\phi\text{ is a partial probability preserving
transformation}\}.$$
It's not hard to see that
$$[T]=\{R:R\text{ a partial probability preserving transformation,}\ \&\ (x,Rx)\in\Cal R\ \text{ a.e.}\}.$$
For $R=T^\phi\in[T]$, write $\phi^{(R)}:=\phi$.
Let
$$[T]_+ =\{R\in [T]:\ \phi^{(R)}\ge 1\text{ a.e.}\}.$$
Recall from \cite{Halm1}:
\proclaim{E.\ Hopf's Equivalence lemma}
If $T$ is an ergodic measure preserving transformation of $(X,\B,m)$ and $A,\ B\in\B$ with
$m(A)=m(B)$, then
$$\exists\ R\in [T]_+\ \text{ such that }\dom(R)=A,\ \Im(R)=B.$$
\endproclaim
We'll also need a quantitative version of this lemma when $A=B$.
\proclaim{Lemma 3.1}
\f Suppose that $T$ is an ergodic probability
preserving transformation of $(X,\B,m)$,
$A\in\B_+$, and $c,\ \e>0,$
then $\forall\ p,q\in\N$ large enough, $\exists\ R\in [T]_+$ such that
$$\dom(R),\ \Im(R)\subset A,\ m(A\setminus\dom(R))<\e,\ \text{ and }\phi^{(R)}=cpq(1\pm\e).$$
\endproclaim
The proof of lemma 3.1 will be given at the end of this section.

\

Let $\Cal R$ be the equivalence relation generated by $T$.
An {\it orbit cocycle } is a measurable function $\tilde\v:\Cal R\to
 G$ such that if $(x,y),\ (y,z)\in\ \Cal R$, then
$$\tilde\v(x,z)=\tilde\v(x,y)+\tilde\v(y,z).$$
\par Let $\v:X\to G$ be measurable, and
let $\v_n\ (n\in\Z)$ denote the cocycle generated
by $\v$ under $T$. The orbit cocycle $\tilde\v:\Cal R\to G$ corresponding to
$\v$ is defined by
$$\tilde\v(x,T^nx)=\v_n(x).$$
For $R\in[T]$, the function $\v_R:\dom(R)\to G$ is defined by
$$\v_R(x)=\tilde\v(x,Rx).$$
Clearly $\v(R\circ S,x)=\v(S,x)+\v(R,Sx)$ on
$\dom(R\circ S)=\dom(S)\cap S^{-1}\dom(R)$.

\

\subheading{\it Definition} Let  $\a$ be a
measurable partition of $X,\ U$ a subset of $G$, and $\e>0$.
We say that the measurable cocycle $\v:X\to\G$
satisfies EVC$_T(U,\e,\a)$ if
\f for $\e$-almost every $a\in \a$, $\exists\ R=R_a\in [T]_+$
such that
$$\dom(R),\ \Im(R)\subset a,\ \v_R\in U\text{ on }\dom(R_a),\ m(\dom(R)))>(1-\e)m(a).$$
\subheading{\it Definition}
We say that the partitions $\{\a_k:k\ge 1\}$
{\it approximately generate} $\B$ if
$$\forall\ B\in\B(X),\ \e>0\ \exists\ k_0\ge 1\ni\ \forall\ k\ge k_0,
\ \exists\ A_k\in\A(\a_k)\ \ni\
m(B\D A_k)<\e.$$
Here $\A(\a)$ denotes the algebra generated by $\a$.
It is not hard to see that
the partitions $\{\a_k:k\ge 1\}$
approximately generate $\B$, if and only if
$E(1_B|\A(\a_k))\to 1_B$ in probability $\forall\ B\in\B$, and
in this case,
$$\forall\ \e>0,\ B\in\B,\ \exists\ k_0
\text{ such that }
\sum_{a\in\a_k,\ 1-m(B|a)\le\e}m(a)\ge (1- \e)m(B)\ \forall\ k\ge k_0.$$
\proclaim{Proposition 3.1}
Suppose that the partitions $\{\a_k:k\ge 1\}$
approximately generate $\B$,
and let $\e_k\downarrow 0,\ \g\in\G$, and $U_k\subset G$ satisfy
$U_n\downarrow \{\g\}$, and $\text{ diam}\, U_n\downarrow 0$.
\par If $\v$ satisfies  EVC$_T(U_k,\e_k,\a_k$)
$\forall\ k\ge 1$, then
$$\g\in E(\v).$$
\endproclaim
\demo{Proof}
Suppose that $B\in\B_+$, and that $V\subset G$ is an open neighbourhood of $\g$. We'll show that
$$\exists\ n\ge 1\ \ni\ m(B\cap T^{-n}B\cap [\v_n\in V])>0.$$
Evidently, $V\supset U_k$ for all $k$ sufficiently large.
It follows from the definitions,
that $\forall\ k$ sufficiently large, $\exists\ a\in\a_k$
such that
$$m(a\setminus B)<0.1 m(a),$$
and $\exists\ R=R_a\in [T]_+$
such that
$$\dom(R),\ \Im(R)\subset a,\ \ \v_R\in U_k\text{ on }\dom(R),\text{ and }m(a\setminus \dom(R))<0.1m(a).$$
It follows that
$$m(B\setminus \dom(R))<0.2m(a).$$
Let $R=T^{\phi}$, where $\phi:\dom(R)\to\Z$. We have that
$$\align &
\sum_{n\in\Z}m(B\cap [\phi=n]\cap T^{-n}B\cap [\v_n\in U_k]) \\ &\ge
m(B\cap\dom(R)\cap R^{-1}(B\cap\Im(R))\cap [\v_R\in U_k])\\ &\ge
0.6m(a),\endalign$$
whence $\exists\ n\in\Z$ such that
$$m(B\cap T^{-n}B\cap [\v_n\in V])\ge
 m(B\cap [\phi=n]\cap T^{-n}B\cap [\v_n\in U_k])>0.$$
\hfill\qed\enddemo
\proclaim{Corollary 3.2} Suppose that the partitions $\{\a_k:k\ge 1\}$
approximately generate $\B$, let $\{U_k:k\ge 1\}$ be a basis of neighbourhoods for the
topology of $G$,
and let $\e_k\downarrow 0$.
\par If $\v$ satisfies  EVC$_T(U_k,\e_k,\a_k$)
$\forall\ k\ge 1$, then $T_\v$ is ergodic.
\endproclaim
This sufficient condition for ergodicity is actually necessary.
\proclaim{Proposition 3.3 } If $T_\v$ is ergodic, then $\forall\
A\in\B_+\ U\ne\emptyset$ open in $G,\ \exists\ R\in [T]_+$ such
that $$\dom(R)=\Im(R)=A,\ \&\ \v_R\in U\text{ a.e. on }A,$$ and
hence, $\v$ satisfies EVC$_T(U,\e,\a)$ for any measurable
partition $\a$ of $X,\ U$ open in $G,\ \e>0$.
\endproclaim
\demo{Proof} This follows from the ergodictity of $T_\v$.
Let $U$ be open in $G$. Choose $g\in U$, then $V:=U-g$ is a neighbourhood of $0\in G$.
Choose $W$ open in $G$ such that $W+W\subset V$.
By ergodicity of $T_\v$, for every $A,B\in\B_+$,
$\exists\ n\in\N$ such that $\mu((A\x W)\cap T_\v^{-n}(B\x (W+g)))>0$,
whence $m(A\cap T^{-n}B\cap[\v_n\in U])>0$.
The proposition follows from this via a standard exhaustion argument.
\qed\enddemo
We'll need a finite version of EVC more suited to sequential constructions.
\subheading{\it Definition} Let  $\a$ be a
measurable partition of $X,\ U$ open in $G,\  \e>0$, and $N\ge 1$.
We say that the measurable cocycle $\v:X\to G$
satisfies EVC$^T(U,\e,\a,N$) if:
\f for $\e$-almost every $a\in \a$, $\exists\ R=R_a\in [T]_+$ with
$\phi^{(R)}\le N$ such that
$$\dom(R),\ \Im(R)\subset a,\ \v_R\in U\text{ on }\dom(R),\text{ and }
m(a\setminus\dom(R))<\e m(a).$$
\proclaim{Proposition 3.4}  Let  $\a$ be a
measurable partition of $X,\ U$ open in $G,\  \e>0$.
The measurable cocycle $\v:X\to\G$
satisfies EVC$_T(U,\e,\a$) iff it satisfies EVC$^T(U,\e,\a,N$) for some
$N\ge 1$.
\endproclaim
The next lemma shows that addition  of a sufficiently small
cocycle does not affect EVC$^T$ conditions too much.
\proclaim{Lemma 3.5}
Let $\a$ be a partition, $\e,\d>0,\ N\in\N$,
$V\subset G$,
and $\phi:X\to G$ be a cocycle
satisfying EVC$^T(U,\e,\a,N$) where $U\subset G$.
\par If $\v:X\to\G$ is  measurable, and
$$m([\v\notin V])<{\d^2\over N},$$
then $\phi+\v$ satisfies
EVC($U+V,\e+\d,\a,N$).
\endproclaim
\demo{Proof}
Let $B=[\v\circ T^j\in V\ \forall\ 0\le j\le N-1]$,
then
since
$$\v_n\in V\text{  on  }B\ \forall\ 1\le n\le N,$$
$$\v_R\in V\text{ on }B\cap\dom(R)\ \forall\ R\in [T]_+\text{ with
}\phi^{(R)}\le N.$$
Let $\a_1$ consist of those $a\in\a$ such that $\exists\ R=R_a\in [T]_+$ with
$\phi^{(R)}\le N$ such that
$$\dom(R),\ \Im(R)\subset a,\ \v_R\in U\text{ on }\dom(R),\text{ and }
m(a\setminus\dom(R))<\e m(a).$$
We have that
$$m(\bigcup_{a\in\a_1}a)>1-\e.$$

Let $\a_2$ consist of those $a\in\a$ for which
$$m(B\cap a)>(1-\d)m(a).$$
It follows from Chebyshev's inequality that
$$m(\bigcup_{a\in\a_2}a)>1-{m(B)\over\d}>1-\d.$$
Therefore
$$m(\bigcup_{a\in\a_1\cap\a_2}a)>1-\e-\d.$$
If $a\in\a_1\cap\a_2$, and $R'=R'_a:=(R_a,\dom(R_a)\cap B)\in [T]_+$,
then:
$$\dom(R'),\ \Im(R')\subset a,\ (\phi+\v)_{R'}\in U+V\text{ on }\dom(R'),\text{ and }
m(a\setminus\dom(R'))<(\e+\d) m(a).$$
\hfill\qed\enddemo
Our main result in this section is a sufficient condition for a group element
to be an essential value of a sum of coboundaries.
\proclaim{Theorem 3.6}
\par Suppose that $g\in G$, the partitions $\{\a_j\}$ approximately
generate $\B$;
\f $N_k\in\N,\ N_k\uparrow\infty$, and $\e_k>0,\ \sum_{k\ge 1}\e_k<\infty.$
\par If for $k\in\N$, $f_k:X\to G$ is measurable,
$$\sum_{j=1}^k(f_j\circ T-f_j)
\text{ satisfies  EVC}^T(N(g,\e_k),\e_k,\a_k,N_k),$$
and
$$m([|f_k\circ T-f_k|\ge {\e_{k-1}\over N_{k-1}}])\le{\e_{k-1}^2\over
N_{k-1}},$$
then
$$\sum_{k=1}^\infty|f_k\circ T-f_k|<\infty\text{   a.e., and  }
g\in E\(\sum_{k=1}^\infty(f_k\circ T-f_k)\).$$
\endproclaim
\demo{Proof}
By the Borel Cantelli lemma,
$\sum_{k=1}^\infty|f_k\circ T-f_k|<\infty$ a.e.. Write
$$\phi:=\sum_{k=1}^\infty(f_k\circ T-f_k),\ \tilde\phi_k=
\sum_{j=1}^k(f_j\circ T-f_j),\ \hat\phi_k=\sum_{j=k+1}^\infty
(f_j\circ T-f_j).$$

\

Since
$$\phi=\tilde\phi_k+\hat\phi_k\ \forall\ k\ge 1,$$
$\tilde\phi_k$ satisfies EVC$^T(N(g,\e_k),\e_k,\a_k,N_k),$
and
$$\align m([|\hat\phi_k|\ge {1\over N_k}\sum_{j=k+1}^\infty\e_j]) &
\le \sum_{j=k+1}^\infty m([f_j\circ T-f_j|\ge {\e_j\over N_k}])\\ &\le
\sum_{j=k+1}^\infty m([f_j\circ T-f_j|\ge {\e_j\over N_{j-1}}\\ & <
\sum_{j=k+1}^\infty {\e_{j-1}^2\over N_{j-1}}\\ &\le
{1\over N_k}\sum_{j=k}^\infty\e_k^2,\endalign$$
it follows from lemma 3.5 that $\phi$ satisfies
\f EVC$^T(N(g,\sum_{j=k}^\infty\e_j),2\sqrt{\sum_{j=k}^\infty\e_k^2},\a_k,N_k)$.
\hfill\qed\enddemo
As promised above, we conclude this section with the
\demo{Proof of lemma 3.1}
Let
$$A_n=\biggl[\bigg|{1\over n}\sum_{k=0}^{n-1}1_A\circ T^k-m(A)\bigg|<\e m(A)\biggr].$$
\f By Birkhoff's ergodic theorem, $\exists\ p_0\in\N$ such that
$m(A_p^c)<{\e^4\over 2}\ \forall\ p\ge p_0.$
Fix $p\ge p_0$. Now fix $q\ge {p\over c\e}:=q_0$. Set
$$B=A_p\cap T^{-[cq]p}A_p.$$
Evidently $m(B)>1-\e^2.$
\f By Birkhoff's ergodic theorem
$\exists\ N_0\in\N$ such that
$$m(C_n^c)<{\e^2\over 2p}\ \ \forall\ n\ge N_0$$
where
$$C_n=\biggl[{1\over n}\sum_{k=0}^{n-1}
1_B\circ T^{pk}\ge E(1_B|\Cal I_{T^p})-\e^2\biggr].$$
\f Let $N>{pq\over\e}\vee pN_0$.
By Rokhlin's theorem, $\exists\ F\in\B$ such that
$$\{T^jF:0\le j\le N-1\}
\text{ are disjoint, and }m\(X\setminus\bigcup_{j=0}^{N-1}T^jF\)<{\e\over p}.$$
Note that since $E(1_B|\Cal I_{T^p})$ is $T^p$-invariant, we have
$$\align {N\over p}\sum_{k=0}^{p-1}\int_{T^kF}E(1_{B^c}|\Cal I_{T^p})dm
& \le\int_XE(1_{B^c}|\Cal I_{T^p})dm\\ &=m(B^c)<\e^2,\endalign$$
whence $\exists\ 0\le k\le p-1$ such that
$$\int_{T^kF}E(1_{B^c}|\Cal I_{T^p})dm<\e^2m(F).$$
There is no loss of generality in assuming $k=0$ as this merely involves
taking $T^kF$ as the base for a slightly shorter Rokhlin tower, and
adding $\bigcup_{j=0}^{k-1}T^jF$ to the "error set".

Set
$$X_0=\bigcup_{j=0}^{N-pq}T^jF,\ \text{ and  }J=X_0\cap \bigcup_{j\ge 0,\ jp\le N}T^{jp}F,$$
 then $m(J)>{1\over 2p}$ so
$$m(C^c_n\cap J)\le \e^2m(J)\ \ \forall\ n\ge N_0.$$
\f For $y\in J$,
set $\k(y)=\#\{0\le j\le p-1:T^jy\in A\}$ and write
$$\{T^jy:0\le j\le p-1,\ T^jy\in A\}=\{T^{j_i(y)}y:1\le
i\le \k(y)\}$$
in case $\k(y)\ge 1$, where $j_i(y)<j_{i+1}(y)$.
Note that
$$\k=pm(A)(1\pm\e)\text{  on  }J\cap A_p.$$
To estimate $m(J\cap B)$:
$$\align \sum_{0\le j\le {N\over p}:\ m(C^c_{{N\over p}}|T^{jp}F)\ge\e}m(T^{jp}F)
& \le \sum_{0\le j\le {N\over p}}{m(C^c_{{N\over p}}\cap T^{jp}F)\over\e}\\ & =
{m(C^c_{{N\over p}}\cap J)\over\e}\\ &\le {m(C^c_{{N\over p}})\over\e}\le {\e\over 2p}\le \e m(J),
\endalign$$
whence, $\exists\ i\le {\e N\over p}$ such that
$$m(C_{{N\over p}}\cap T^{ip}F)=m(T^{-ip}C_{{N\over p}}\cap F)\ge (1-\e)m(F).$$
For $y\in T^{-ip}C_{{N\over p}}\cap F$,
$$\align \#\{0\le j\le {N\over p}:T^{jp}y\in B\} & \ge
\#\{0\le j\le {N\over p}:T^{(i+j)p}y\in B\}-\e{N\over p}\\ &
\ge {N\over p}(E(1_B|\Cal I_{T^p})-2\e).\endalign$$
Therefore,
$$\align m(J\cap B) &=\sum_{k=0}^{{N\over p}-1}m(T^{jp}F\cap B)\\
&= \int_F\(\sum_{k=0}^{{N\over p}-1}1_{B}\circ
T^{jp}\)dm\\ &\ge
{N\over p}\int_{T^{-ip}C_{{N\over p}}\cap F}(E(1_B|\Cal I_{T^p})-2\e)dm\\ &
\ge
{N\over p}\int_F(E(1_B|\Cal I_{T^p})-3\e)dm\\ &
\ge (1-4\e)m(F){N\over p} = (1-4\e)m(J).\endalign$$
\f For $x\in\bigcup_{j=0}^{p-1}T^jJ$,
let $j(x)$ be such that $T^{-j(x)}x\in J$, and let
$y(x)=T^{-j(x)}x$.
\f Define
$\psi:A\cap\bigcup_{j=0}^{p-1}T^jJ\to \{1,2,\dots p\}$ by
$$\psi(x)=\sum_{k=0}^{j(x)}1_A(T^{-k}x)=\sum_{k=0}^{j(x)}1_A(T^ky(x)).$$
Note that
$$x=T^{j_{\psi(x)}(y(x))}y(x).$$
\par Now define $D\subset A\cap X_0$ by
$$D\cap\bigcup_{j=0}^{p-1}T^jJ_0=
\{x\in A\cap J_0:\psi(x)\le
\k(y(T^{[cq]p}x))\},$$
and define $\phi:D\to\N$ by
$$\phi(x)=[cq]p+j_{\psi(x)}(y(T^{[cq]p}x))\ \ \ x\in
D\cap\bigcup_{j=0}^{p-1}T^jJ.$$
\par We claim that if $R\in [T]_+$ is defined by $\dom(R)=D$ and
$\phi^{(R)}=\phi$, then $\phi$ is as advertised.
To see this,
check that $\k\ge (1-\e)m(A)p$ on $J\cap B$, whence
$$m(D)\ge m(J_0\cap B)(1-\e)m(A)p\ge (1-6\e)m(J)pm(A)\ge (1-7\e)m(A).$$
\hfill\qed\enddemo
\heading
\S4  proof of theorems 1 and 2
\endheading
In this section, we prove theorems 1 and 2. The proofs are sequential
using theorem 3.6. The inductive steps are lemmas 4.1 and 4.2.
Their proofs use the Rohlin lemmas for Abelian group actions
of Katznelson and Weiss \cite{K-W}, and Lind \cite{Lin} respectively
(see also \cite{O-W} for a general Rohlin lemma for amenable
group actions implying these).

\

\par Let $G$ be a locally compact, second countable Abelian group
with invariant metric $d$,
and let $T$ be an ergodic probability preserving transformation of the standard probability space $(X,\B,m)$.
\subheading{\it Definition}
\par A measurable function $f:X\to G$ is called a {\it $T$-coboundary}
if $f=h-h\circ T$ for some  measurable function $h:X\to G$.
\par Measurable functions $f,g:X\to G$ are said to be {\it $T$-cohomologous},
written $f\overset{T}\to\sim g$,  if $f-g$ is a $T$-coboundary.

\

Let $\v:X\to G$ be measurable. Suppose $S\in\L_\v(T)$, and $w\in\EE(G)$,
then
$$w=\pi_\v(S)\ \ \Leftrightarrow\ \ \v\circ S\overset{T}\to\sim w\circ\v.$$
We prove the following version of theorem 1:
\proclaim{Theorem 1'}
\f Suppose that $T$ is an ergodic probability
preserving transformation, $S_1,\dots,S_d\in
C(T)\ (d\le\infty)$ are such that $(T,S_1,\dots,S_d)$ generate a free
$\Z^{d+1}$ action of probability preserving transformations of $X$.
\par If $w_1,\dots,w_d\in\EE(G)$ commute
 (i.e. $w_i\circ w_j=w_j\circ w_i\ \forall\ 1\le i,j\le d$), then there is
a measurable function $\v:X\to G$ such that
$T_\v$ is ergodic, and
$$\v\circ S_i\overset{T}\to\sim w_i\circ\v\ \ (1\le i\le d).$$
\endproclaim

\

\proclaim{Lemma 4.1} Let $\phi:X\to G$ be a $T$-coboundary, let
$S_1,\dots,S_d$ be probability preserving transformations generating a free $Z^{d+1}$ action together with
$T$, and let $w_1,\dots,w_d\in\EE(G),\ w_i\circ w_j= w_j\circ w_i$.
If $\a$ is a finite,
measurable partition of $X$, and $\e>0,$
then there is a measurable function $f:X\to G$ such that
$$m([f\circ T-f\neq 0])<\e,\tag1$$
$$m([f\circ S_j\neq w_j\circ f])<\e,\ (1\le j\le d)\tag2$$
and such that
$$\phi+f-f\circ T\text{ satisfies EVC}_T(N(\g,\e),\e,\a).\tag3$$
\endproclaim
\demo{Proof} Write $\phi=H-H\circ T$. Possibly refining $\a$, we may
assume that for ${\e\over 2}$-a.e. $a\in\a$, the oscillation of $H$ on
$a$ is less than ${\e\over 2}$.
\par For $\u i=(i_1,\dots,i_d)\in\Z_+^d$, we'll write
$$S_{\u i}:=S_1^{i_1}\circ\cdots\circ S_d^{i_d},\ \& \
w_{\u i}:=w_1^{i_1}\circ\cdots\circ w_d^{i_d}.$$
Then
$$S_{\u i+\u j}=S_{\u i}\circ S_{\u j},\ \& \
w_{\u i+\u j}=w_{\u i}\circ w_{\u j}$$
since $S_i\circ S_j= S_j\circ S_i$ and $w_i\circ w_j= w_j\circ w_i$.
\par Fix $k>{10\over\e}$.
There is an ergodic cocycle
$\v:X\to G$ such that
$$m([\v\neq 0])<{\e\over 3k^d}.$$
It follows that $w_{\u i}\circ\v\circ S_{-\u i}$ is ergodic for $\u i\ge\u 0$
(as $w_{\u i}$ is surjective, and $S_{-\u i}$ commutes with $T$
 for $\u i\ge\u 0$), whence
$\phi+ w_{\u i}\circ\v\circ S_{-\u i}$ is ergodic for $\u i\ge\u 0$
(as $\phi$ is a coboundary), and so satisfies
EVC$_T(N(\g,{\e\over 4}),{\e\over 4k^{d}},\a)$.
Therefore (by propositions 3.3 and 3.4),
$$\exists\ M\in\N\text{  such that }
\phi+ w_{\u i}\circ\v\circ S_{-\u i}
\text{  satisfies EVC}^T(N(\g,{\e\over 4}),{\e\over 4k^d},\a,M)$$
 for $\u 0\le\u i\le\u k$ where
$\u k:=(\undersetbrace\text{$d$-times}\to{k,\dots,k})$,
and
$\u 0\le \u i< \u k$ means $0\le i_j< k_j\ \forall\ 1\le j\le d$.

\

Now choose $N\ge 1$ such that
$${M\over N}<{\e\eta_\a\over 5}$$
where $\eta_\a:=\min\,\{m(a):a\in\a\}$.
By the Katznelson-Weiss Rohlin lemma (\cite{K-W}),
$\exists\ F\in\B(X)$ such that
$$\{T^jS_{\u i}F:\ 0\le j\le N-1,\  \u 0\le \u i< \u k\} \text{ are disjoint,}$$
and
$$m\(X\setminus\bigcup_{0\le j\le N-1,\  \u 0\le \u i< \u k}
T^jS_{\u i}F\)<{\e\eta_\a\over
6}.$$

\

Let
$$C=\bigcup_{j=0}^{N-1}T^jF,\ \ \widetilde C=\bigcup_{j=0}^{N-M}T^jF,
\ \ \Cal T=\bigcup_{\u 0\le\u i<\u k}S_{\u i}C,
\ \ \widetilde{\Cal T}=\bigcup_{\u 0\le\u i<\u k}S_{\u i}\widetilde C.$$
There is a measurable function $f_0:X\to G$ such that
$$\v=f_0 -f_0\circ T
\text{ on }\Cal T.$$
Set $\v'=f_0 -f_0\circ T$, then $m([\v\ne\v'])<{\e\eta_\a\over 6}$.
\

Now define $f:\Cal T\to G$ by
$$f =\cases &  w_{\u i}\circ f_0\circ S_{-\u i}\ \ \  \text{ on }S_{\u i}C\ \
(\u 0\le\u i\le\u k)\\ & 0\ \ \ \text{else,}\endcases$$
and define
$$\psi=f-f\circ T.$$

\

To establish (1),
$$\align m([\psi\ne 0]) &<m([\psi\ne 0]\cap \widetilde{\Cal T})
+m(X\setminus\widetilde{\Cal T})\\ &\le
k^dm([\v\ne 0]\cap \widetilde C)+m(X\setminus\widetilde{\Cal T})
\\ &<{\e\over 3}+{M\over N}\\ &<\e.\endalign$$

\

Next, to prove (2), suppose that $\u 0\le\u i<\u k,\ 1\le j\le d$ and
$i_j<k-1$. If $x\in S_{\u i}C$, then
$$\align f(S_jx) & =
w_{\u i+\u e_j}\circ f_0\circ S_{-(\u i+\u e_j)}(S_jx)\\ &=
w_j\circ w_{\u i}\circ f_0\circ S_{-\u i}(x)\\ &=
w_j\circ f(x)\endalign$$
whence
$$\align m([f\circ S_j\ne w_j\circ f]) &<m\(\bigcup_{\u 0\le\u i<\u k,\ i_j=k-1}
S_{\u i}C\)+m(X\setminus\Cal T)\\ &<{1\over k}+{\e\eta_\a\over 6}<\e.
\endalign$$

\

To complete the proof, we show (3).
\f We know that $\phi+w_{\u i}\circ \v\circ S_{-\u i}$ satisfies
EVC$^T(N(\g,{\e\over 4}),{\e\over 4k^{d}},\a,M)\ \forall\ \u i$, whence
\f $\phi+w_{\u i}\circ \v'\circ S_{-\u i}$ satisfies
EVC$^T(N(\g,{\e\over 4}),{\e\over 3k^{d}},\a,M)\ \forall\ \u i$.
\par It follows that  for ${\e\over 3}$-a.e. $a\in\a$, and
for each $\u 0\le\u i<\u k$,
\f $\exists\ R_{\u i}=R_{a,\u i}\in [R]_+$ such that $\dom(R_{\u i}),\ \Im(R_{\u i})\subset a$,
\f $m(a\setminus\dom(R_{\u i}))<{\e\over 3k^d}m(a)$, and
\f $(\phi+w_{\u i}\circ \v\circ S_{-\u i})_{R_{\u i}}\in N(\g,{\e\over 4})$ on $\dom(R_{\u i})$.
\par Define $R=R_a\in [T]_+$ by
$$\dom(R)=\bigcup_{\u 0\le\u i<\u k}\dom(R_{\u i})\cap S_{\u i}\widetilde C,$$
and
$$R=R_{\u i}\text{ on }S_{\u i}\widetilde C,\ \ (\u 0\le\u i<\u k).$$
For $x\in\dom(R),\ \exists\ \u i=\u i(x)$ such that $x\in\dom(R_{\u i})\cap S_{\u i}\widetilde C,$
and  we have that
$$\align (\phi+\psi)_R(x) & =(\phi+w_{\u i}\circ \v'\circ S_{-\u i})_{R_{\u i}}(x)\\
&\in N(\g,{\e\over 4}).\endalign$$
Lastly,
$$\align m(a\setminus\dom(R)) & =\sum_{0\le\u i<\u k}m((a\setminus\dom(R))\cap S_{\u i}\widetilde C)
+m(\Cal T\setminus\widetilde{\Cal T)}+m(X\setminus\Cal T)\\
&<\sum_{0\le\u i<\u k}m(a\cap S_{\u i}\widetilde C\setminus\dom(R_{\u i}))+{M\over N}+m(X\setminus\Cal T)\\
&\le\sum_{0\le\u i<\u k}m(a\setminus\dom(R_{\u i}))+{\e\over 5}\eta_\a+{\e\over 6}\eta_\a\\ &\le \e m(a).\endalign$$
\hfill\qed\enddemo
\demo{Proof of theorem 1' in case $d$ finite}
\par Choose a countable, dense subset $\G$ of $G$.
Let $(\g_1,\g_2,\dots)\in\G^\N$ satisfy
$$\{\g_k:k\ge 1\}=\G,\ \&\ \forall\ \g\in\G,\ \g_k=\g\text{ for
infinitely many }k,$$
let the partitions $\{\a_j\}$ approximately
generate $\B$, and let $\e_k=2^{-k^2}$.

Construct (sequentially) using lemma 4.1, a sequence of
coboundaries
$$\phi_k=f_k-f_k\circ T$$
such that
$$m([f_k\circ S_j\ne E\circ f_k])\le\e_k\ (1\le j\le d),$$
$\tilde\phi_k:=\sum_{j=1}^k\phi_j$ satisfies EVC$^T(N(\g_k,\e_k),\e_k,\a_k,N_k$) where
$N_k\in\N,\ N_k\uparrow$, and
$$m([\phi_k\neq 0])\le {\e_k\over N_{k-1}}.$$
Clearly
$\phi:=\sum_{k=1}^\infty\phi_k$ converges a.e.. Also
$$\psi_j:=\sum_{k=1}^\infty (f_k\circ S_j- w_j\circ f_k)\ (1\le j\le d)$$
 converges a.e., whence
$$\phi\circ S_j- w_j\circ\phi=\psi_j-\psi_j\circ T\ (1\le j\le d).$$
Theorem 3.6 now shows that $\G\subset E(\phi)$, and the
ergodicity of $\phi$ is established.
\hfill\qed\enddemo
We prove the following version of theorem 2.
\proclaim{Theorem 2'} Suppose that $T$ is an ergodic probability
preserving transformation, $\{S_t:t\in\R\}
\subset C(T)$ are such that $T$ and $\{S_t:t\in\R\}$ generate a free
$\Z\x\R$ action of probability preserving transformations of $X$.
\par There is
a measurable function $\v:X\to\R$ such that
\f $T_\v$ is ergodic; and
\f $\exists\ g:\R\x X\to\R$ measurable (with respect to $m_\R\x m$)
such that
$$\v\circ S_t(x)-e^t\v(x)=g(t,Tx)-g(t,x),\tag1$$
and
$$g(t+u,x)=g(t,S_ux)+e^tg(u,x).\tag2$$
\endproclaim
If $Q_t(x,y):=(S_tx,e^ty+g(t,x)$, then
(1) implies that $Q_t\in C(T_\v)\ \forall\ t\in\R$, and
by (2), $\{Q_t:t\in\R\}$ is a flow, whence
$T_\v$ is a Maharam transformation.
\proclaim{Lemma 4.2} Let $\phi:X\to\R$ be a $T$-coboundary, let
$\{S_t:t\in\R\}$ be probability preserving transformations generating a free $\Z\x\R$ action together with
$T$.
\par If $\a$ is a finite,
measurable partition of $X$, $\e>0$, and $J\subset\R_+$ is an open interval,
then there is a measurable function $f:X\to\R$ such that
$$m([|f\circ T-f|\ge\e])<\e,\tag1$$
$$m([f\circ S_t\neq e^tf])<\e,\ (0\le t\le 1)\tag2$$
and such that
$$\phi+f-f\circ T\text{ satisfies EVC}_T(J,\e,\a).\tag3$$
\endproclaim
\demo{Proof of lemma 4.2}
Write $J=((1-\d)b,(1+\d)b)$ where $b,\ \d>0$.
We'll sometimes use the notation $x=(1\pm\d)b$ which means $x\in J$.
\par Write $\phi=\psi\circ T-\psi$ where $\psi:X\to\R$ is measurable.
\par Choose a refinement $\a_1$ of $\a$ with the property that
$$\forall\ a\in\a_1,\ \exists\ y_a\in\R\ \ni\ |\psi-y_a|<{b\d\over
2}\text{ a.e. on }a,$$
and set $\eta_\a:=\min\,\{m(a):a\in\a\}$.
\par Fix $K={10\over\e}$, and $0=t_0<t_1<\dots<t_M=K$ such that
$e^{t_{i+1}}<(1+{\d\over 3})e^{t_i}$.

\

By lemma 3.1, $\exists\ p,q\in\N$ such that ${be^K\over pq}<\e$, and
$$\forall\ a\in\a_1,\ \ 0\le k\le M-1,\ \exists\ R_{a,k}\in [T]_+$$
such that
$$\dom(R_{a,k}),\ \Im(R_{a,k})\subset a,\ m(a\setminus\dom(R_{a,k}))
<{\e\over 7M}m(a),
\text{ and }\phi^{(R_{a,k})}=e^{-t_k}pq(1\pm {\d\over 9}).$$

\

Now choose $N\ge 1$ such that
$${e^Kpq\over N}<{\e\eta_\a\over 5}.$$
By the Rokhlin theorem for continuous groups
(\cite{Lin}, \cite{O-W})
$$\exists\ F\in\B(X)\text{ such that }T^kS_tF\text{ are disjoint for }0\le k\le N,\ 0\le t\le K,$$
and
$$m\(X\setminus\bigcup_{0\le k\le N-1,\ 0\le t\le K}
T^kS_tF\)<{\e\eta_\a\over 6}.$$

\

Let
$$C=\bigcup_{j=0}^{N-1}T^jF,\ \ \widetilde C=\bigcup_{j=0}^{N-2}T^jF,
\ \ \Cal T=\bigcup_{0\le t\le K}S_tC,\ \
\widetilde{\Cal T}=\bigcup_{0\le t\le K}S_t\widetilde C.$$

There is a measurable function $f:\Cal T\to\R$ such that
$$f\circ T -f= {b\over pq}e^t
\text{ on }S_t\widetilde C.$$
Complete the definition of $f:X\to\R$ by setting $f=0$ on $\Cal T^c$.

\

It is immediate from this construction that $f$ satisfies (1) and (2). We establish (3) by showing
that $f\circ T-f$ satisfies EVC$_T(J,\e,\a_1)$.
Let
$$\widehat C=\bigcup_{j=0}^{N-pq}T^jF,\ \ \widehat{\Cal T}=\bigcup_{0\le t\le K}S_t\widehat C.$$
Let, for $0\le k\le M-1$,
$$\widehat{\Cal T}_k=\bigcup_{t_k\le t<t_{k+1}}S_t\widehat C.$$
Fix $a\in\a_1$, and define $R'_a\in [T]_+$
by
$$R'_a=R_{a,k}\text{ on }\dom(R_{a,k})\cap\widehat{\Cal T}_k.$$
It follows that $\dom(R'_a),\ \Im(R'_a)\subset a$,
$$\align m(a\setminus\dom(R'_a)) & =
\sum_{k=0}^{M-1}m(\widehat{\Cal T}_k\cap [a\setminus\dom(R_{a,k})])\\ &\le
\sum_{k=0}^{M-1}m(a\setminus\dom(R_{a,k}))\\ &\le
{\e\over 7}m(a);\endalign$$
and, on $\dom(R'_a)\cap\widehat{\Cal T}_k$:
$$|\psi\circ R'_a-\psi|<{b\d\over 2},$$
whence, on $S_t\widetilde C,\ t\in [t_k,t_{k+1}]$,
$$\align \v_{R'_a} &={e^tb\over pq}\phi^{(R'_a)}\pm {b\d\over 2}\\ &=
e^{t-t_k}b(1\pm{\d\over 9})\pm {b\d\over 2}\\ &=
b(1\pm{\d\over 9})(1\pm {\d\over 3})(1\pm {\d\over 2})
\in J.\endalign$$
\hfill\qed\enddemo
\demo{Proof of theorem 2'}
Fix $(g_1,g_2,\dots)=(1,\sqrt 2,1,\sqrt 2,\dots)$.
\f Construct using lemma 4.2, a sequence of
coboundaries
$$f_k\circ T-f_k$$
such that
$$m([f_k\circ S_t\ne e^tf_k])\le{1\over 2^k}\ \ \ (0\le t\le 1),$$
$$\phi_k:=\sum_{j=1}^k(f_j\circ T-f_j\text{ satisfies EVC}^T
((\g_k-{1\over 2^k},\g_k+{1\over 2^k}),\e_k,\a_k,N_k)$$
where
$N_k\in\N,\ N_k\uparrow$, and
$$m([|f_k\circ T-f_k|\geq {1\over 2^kN_{k-1}}])\le {1\over 2^kN_{k-1}}.$$
The ergodicity of
$$\sum_{k=1}^\infty(f_k\circ T-f_k)$$
follows from
$$1,\ \sqrt 2\ \in\ E\(\sum_{k=1}^\infty(f_k\circ T-f_k)\)$$
which follows from theorem 3.6
\hfill\qed\enddemo
\heading\S5 Maharam transformations\endheading
In this section, we give conditions for a
conservative, ergodic, measure preserving transformation to be isomorphic to
a Maharam transformation. The first proposition shows that the transformations constructed
in theorem 2 are Maharam transformations, and the second (a converse of theorem 2 for Kronecker
transformations) will be used to
construct completely squashable, ergodic $\R$-extensions which are not
isomorphic to any Maharam transformation.
\proclaim{Proposition 5.1}
A conservative, ergodic, measure preserving transformation $T$ of
the standard, non atomic, $\s$-finite  measure space $(X,\B,m)$ is isomorphic to
a Maharam transformation if, and only if there is a flow $\{Q_t:t\in\R\}\subset C(T)$
such that $D(Q_t)=e^t\ \forall\ t\in\R$.
\endproclaim
\demo{Proof}
\f Suppose first that $T$ is a Maharam transformation, i.e.
$T:X=\Om\x\R\to X$ is defined by
$$T(x,y)=(Rx,y-\log{dp\circ R\over d\,p})$$
and preserves the measure $dm(x,y):=dp(x)e^ydy$, where
$R$ is a non-singular conservative, ergodic transformation of
the standard probability space $(\Om,\A,p)$.
Set $Q_t(x,y)=(x,y+t)$, then $\{Q_t:t\in\R\}\subset C(T)$
is a flow, and $D(Q_t)=e^t$.
\par Conversely, suppose that there is a flow $\{Q_t:t\in\R\}\subset C(T)$
such that $D(Q_t)=e^t\ \forall\ t\in\R$. The flow $\{Q_t:t\in\R\}$ is
dissipative on $X$. It is well known that up to measure theoretic
isomorphism, $X=\Om\x\R$
where $\Om$ is some probability space,
$Q_t(x,y)=(x,y+t)$, and $dm(x,y)=e^ydp(x)dy$ where $p$ is the probability
on $\Om$.
\par Since $\{Q_t:t\in\R\}\subset C(T)$, $\exists$ a non-singular transformation
$R:\Om\to\Om$ such that
$$T(x,y)=(Rx,Y(x,y)).$$
A calculation shows that indeed
$Y(x,y)=y-\log R'(x)$ where $R'={d\l\circ R\over d\,\l}$,
i.e. $T$ is the Maharam transformation of $R$. The ergodicity of $T$ implies that
$\Om$ is non-atomic, and hence standard.
\hfill\qed\enddemo
\remark{Remarks}
\f 1) By proposition 5.1, the skew products constructed in theorem 2
are isomorphic to Maharam transformations.
\f 2) Let $T$ be a Bernoulli transformation. We claim that there is a
$\Z$-extension of $T$ which is isomorphic to a Maharam transformation.
\par Indeed, by theorem 2 and the above remark,
there is such an $\R$-extension of $T$. By proposition 2.5, this
$\R$-extension of $T$ is isomorphic to a $\Z$-extension of $T$.
\endremark
\proclaim{Proposition 5.2}
\f Let $T$ be a Kronecker transformation of the compact, metric, Abelian group $X$.
\par If there is an ergodic $\R$-extension of $T$ which is isomorphic to some
Maharam transformation, then there is a  continuous, injective group
homomorphism $\R\to X$.
\endproclaim
\demo{Proof}
Let $T_\v$ be an ergodic $\R$-extension of $T$ which is isomorphic to some
Maharam transformation.
By proposition 5.1, there is a flow $\{Q_t:t\in\R\}\subset C(T_\v)$
such that $D(Q_t)=e^t\ \forall\ t\in\R$. It follows from \cite{A-L-M-N} that
$Q_t$ has form $(*)\ \forall\ t\in\R$, i.e.
$$Q_t(x,y)=(S_tx,e^ty+g_t(x)).$$
Clearly, the map $t\mapsto S_t$ is a measurable homomorphism from $\R\to
C(T)\cong X$,
whence by Banach's theorem,
continuous. To see that  $t\mapsto S_t$ is injective,
suppose otherwise, that $S_a=$Id for some $a\ne 0$.
Then $Q_a(x,y)=(x,e^ay+g_a(x))$, whence
$$\v=G-G\circ T\text{ where }G={g_a\over e^a-1}$$
contradicting ergodicity of $T_\v$.
\hfill\qed\enddemo
\heading\S6 completely squashable $\R$-extensions of odometers\endheading
For $a_n\in\N,\ (n\in\N)$, set
$$X:=\prod_{n=1}^\infty\{0,\dots,a_n-1\}$$
equipped with the addition
$$(x+x')_n= x_n+x'_n+\e_n\ \ \mod a_n$$
where
$$\e_1=0,\ \&\ \e_{n+1}=\cases & 0\ \ \ x_n+x'_n+\e_n<a_n\\
& 1\ \ \ x_n+x'_n+\e_n\ge a_n.\endcases$$
Clearly, $X$ equipped with the product discrete topology,
is a compact Abelian topological group, with Haar measure
$$m=\prod_{n=1}^\infty\,({1\over a_n},\dots,{1\over a_n}).$$
Also if $\tau=(1,0,\dots)$ then $X=\o{\{n\tau\}}_{n\in\Z}$
whence $x\mapsto Tx(:=x+\tau)$ is ergodic.

\

Set
$q_1=1,\ q_{n+1}=\prod_{k=1}^na_k,$
then
$$(q_n\tau)_k=\cases & 1\ \ \ k=n\\ & 0\ \ \ k\ne n,\endcases$$
whence
$$T^{q_n}x=(x_1,\dots,x_{n-1},\tilde T_n(x_n,\dots))$$
where
$\tilde T_n:\prod_{k=n}^\infty \{0,\dots,a_k-1\}\to
\prod_{k=n}^\infty \{0,\dots,a_k-1\}$ is defined by
$\tilde T_n(x)= x+\tilde\tau_n$ where
$\tilde\tau_n=(1,0,\dots).$

\

The transformation $T\cong (X,T)$
is called the {\it odometer with digits $\{a_n:n\in\N\}$}.
Let $G$ be a second countable LCA group, and let $(X,T)$ be an odometer.
\par We consider cocycles $\v:X\to G$
of form $$\v(x):=\sum_{n=1}^\infty[\b_n((Tx)_k)-\b_n(x_k)],$$
the sum being a finite sum. Cocycles of this form are called
{\it of product type}. We'll call the functions $\{\b_k:k\in\N\}$
the {\it partial transfer functions} of $\v$.
Clearly if the sum of the partial transfer functions converges, then
indeed the limit is a transfer function for the cocycle, which is a coboundary.
\par We prove
\proclaim{Theorem 6.1}
There is an odometer $(X,T)$, and an
ergodic cocycle $\v:X\to\R$ of product type such that
$T_\v$ is completely squashable, indeed
\f $\forall\ c>0,\ \exists$
a measurable function $\psi_c:X\to G$, and a translation
$S_c:X\to X$ satisfying
$$\v\circ S_c=c\v+\psi_c\circ T-\psi_c.$$
\endproclaim

\

We claim that $T_\v$ is not isomorphic to a Maharam transformation.
\par Otherwise, by proposition 5.2,
there would be a continuous, injective group homomorphism
$\R\to X$, whose existence is prevented by
the disconnectedness of $X$.

We prove ergodicity of $\v$ using rigid essential value conditions.
\proclaim{Proposition 6.2}
Let $\v:X\to G$ be a cocycle, and let $\g\in G$.
If $\forall\ \e>0,\ \ \exists\ \d_k\to 0,$
and a sequence of partitions $\a_k$ which
approximately generate $\B$, such that
$$m\(\bigcup_{a\in\a_k}a\)>1-\d_k,$$
and for every $k\ge 1$, for $\d_k$-a.e. $a\in\a_k,\ \exists\ n=n(a)$
such that
$$m(a\D T^{-n}a)<\d_km(a),\text{ and }m(a\cap [\v_n\in N(\g,\e)])>{m(a)\over 25};$$
then
$$\g\in E(\v).$$
\endproclaim
\demo{Proof} This is a special case of lemma 3.1.\qed
\enddemo
The functions $\b_k$ are defined by means of blocks.
To $\u\g=(\g_1,\g_2,\dots,\g_m)\in\R^m$, associate a
{\it canonical difference block }
$B_{\u\g}=(b(0),b(1),\dots,b(2^m-1))\in\R^{2^m}$
defined  by
$$b\(\sum_{k=1}^m\e_k2^{k-1}\)=
\sum_{k=1}^m\e_k\g_k\ \ \ \ (\e_1,\dots\e_m)\in\{0,1\}^m.$$
It is evident that for $0\le \nu\le 2^m-1,\ 1\le j\le m$ with $\e_j(\nu)=0$,
we have $\nu+2^{j-1}\le 2^m-1$, and
$b(\nu+2^{j-1})-b(\nu)=\g_j.$ It follows that
$\forall\ 1\le j\le m,\ \exists\ 1\le n=n(j)\le 2^m$ such that
$$\#\{1\le\nu\le 2^m-1:\nu+n\le 2^m-1,\ b(\nu+n)-b(\nu)=\g_j\}
\ge {2^m\over 2}.\tag1$$

We'll need some control over the size of
$|b(j)|,\ \ (0\le j\le 2^m-1)$
and to obtain this, we need the
\f {\it balanced canonical difference block }
associated to $\u\g=(\g_1,\g_2,\dots,\g_m)\in\R^m$,
defined by
$$B=(b(0),b(1),,\dots,b(4^m-1))\in\R^{4^m}$$
 where
$$b\(\sum_{k=1}^m\e_k2^{k-1}+\sum_{\ell=1}^m\d_\ell
2^{m+\ell-1}\)=\sum_{k=1}^m(\e_k-\d_k)\g_k,\ \ (\u\e,\u\d\in\{0,1\}^m).$$

\

\par Let  $B$ be the balanced canonical
difference block associated to $(\g_1,\dots,\g_m)\in\R^{m}$. Since
$B$ is also the canonical difference block associated to
\f$(\g_1,\dots,\g_m,-\g_1,\dots,-\g_m)\in\R^{2m}$, we have by (1)
that
$$\#\{1\le\nu\le 4^m-1:\nu+n\le 4^m-1,\ b(\nu+n)-b(\nu)=\g_j\}
\ge {4^m\over 2}.\tag2$$
Also, we claim that
$$\#\{0\le \nu\le 4^m-1:b(\nu)|\ge m^{{3\over 4}}\}\le\max_{1\le j\le m}
|\g_j|^2{4^m\over\sqrt m}.\tag3$$
To see this
$$\align
|\{0\le \nu\le 4^m-1:|b(\nu)|\ge m^{{3\over 4}}\}| &\le
{1\over m^{{3\over 2}}}\sum_{\u\e,\u\d\in\{0,1\}^m}\(\sum_{k=1}^m
(\e_k-\d_k)\g_k\)^2\\
&={4^m\over m^{{3\over 2}}}\sum_{k=1}^m{\g_k^2\over 2}\\ &\le
{4^m\over \sqrt m}\max_{1\le j\le m}|\g_j|^2.\endalign$$
\par We now construct the odometer and cocycle.
We construct our cocycle $\v$ to have $1,{1\over\sqrt 2}\in E(\v)$ thus ensuring ergodicity.
Let $g_{2n}=1$, and $g_{2n+1}={1\over\sqrt 2}$.

\

For $k\ge 1$, choose natural numbers $\nu_k$ and $\mu_k$ satisfying
$$\sum_{k=1}^\infty{1\over\mu_k}<\infty,\tag3$$
$$\sum_{k=1}^\infty{e^{2\mu_k}\over\sqrt{\mu_k\nu_k}}<\infty,\tag4$$
and
$$\sum_{k=1}^\infty{\mu_k\over\nu_k^{{1\over 4}}}<\infty.\tag5$$
For example:
$$\mu_k=k^2,\ \text{ and }\nu_k=k^23^{4k^2}.$$
Set $m_k=\mu_k\nu_k$, and let
$$B_k=(b_k(0),b_k(1),\dots,b_k(4^{m_k}-1))$$
be the balanced canonical difference block associated to
$(\g_k(1),\dots,\g_k(m_k))$ where
$$\g_k(j)=g_ke^{-{j-1\over\nu_k}}.$$
Now let
$$a_k=m_k4^{m_k}$$
and let $(X,T)$ be the odometer with digits $\{a_n:n\in\N\}$.
\par We specify a cocycle $\v:X\to \R$ of product type, defining it's
partial transfer functions
$\b_k:\{0,\dots,a_k-1\}\to \R$ by
$$\b_k(j4^{m_k}+\nu)=e^{{j\over\nu_k}}b_k(\nu)\ \ \ (0\le j\le m_k-1,\ 0\le\nu\le 4^{m_k}-1).$$

\

Note that for $0\le j\le m_k-1$,
$$\align & |\{j4^{m_k}\le\nu<(j+1)4^{m_k}:\b_k(\nu+n(j,k))-\b_k(\nu)=g_k\}|\\ & \ge
|\{0\le\nu<4^{m_k}:b_k(\nu+n(j,k))-\b_k(\nu)=\g_k(j)\}|
\\ & \ge {4^{m_k}\over 2}.\tag6\endalign$$

\

\

To check ergodicity of $\v$, we show, using the essential value
condition that $ 1,{1\over\sqrt 2}\in E(\v)$. For $k\ge 1$, we let
$$\a_k=\{A((u_1,\dots,u_{k-1}),j):0\le u_\nu <a_\nu,\ 1\le\nu\le k-1,\
0\le j\le m_k-2\}$$
where
$$A((u_1,\dots,u_{k-1}),j)=\{x\in X:x_\nu=u_\nu,\ 1\le\nu\le k-1,\ \&\
jN_k\le u_k<(j+1)N_k\}.$$
It follows that
$$m\(\bigcup_{a\in\a_k}a\)=1-{1\over m_k}.$$
Also, if $n=n(j,k)q_k$ where $n(j,k)$ is as in (6), then
$$m(A(\u u,j)\D T^{-n}A(\u u,j))<{1\over m_k}m(A(\u u,j)),$$
and, by (6)
$$m(A(\u u,j)\cap [\v_{n(j,k)}=\g_k])\ge {m(A(\u u,j))\over 2}.$$
The essential value condition now shows that
$$1,{1\over\sqrt 2}\in E(\v).$$

\

To conclude, we show that $\forall\ c\in (1,e),\ \exists
S:X\to X$ such that $\v\circ S-c\v$ is a coboundary.
Fix $c\in (1,e)$ and let
$$r_k=[\nu_k\log c]\le\nu_k.$$
Let
$$S=(r_14^{m_1},r_24^{m_2},\dots).$$
We claim that $\forall\ k\ge 1$,
$$m(\{x\in X:|\b_k(x_k)|\ge m_k^{{3\over 4}}\})\le {e^{2\mu_k}\over\sqrt{m_k}}.$$
To see this
$$\align & m(\{x\in X:|\b_k(x_k)|\ge m_k^{{3\over 4}}\})\\ & =
{1\over a_k}\#\{0\le\nu\le a_k-1:|\b_k(\nu)|\ge m_k^{{3\over 4}}\}\\ &=
{1\over a_k}\sum_{j=0}^{m_k-1}\#\{0\le\nu\le 4^{m_k}-1:|\b_k(j4^{m_k}+\nu)|\ge m_k^{{3\over 4}}\}\\ &=
{1\over a_k}\sum_{j=0}^{m_k-1}\#\{0\le\nu\le 4^{m_k}-1:e^{{j\over\nu_k}}|b_k(\nu)|\ge m_k^{{3\over 4}}\}\\ &\le
{1\over 4^{m_k}}\#\{0\le\nu\le 4^{m_k}-1:e^{\mu_k}|b_k(\nu)|\ge m_k^{{3\over 4}}\}\\ &\le
{e^{2\mu_k}\over\sqrt{m_k}}.\endalign$$
It now follows from (4), and the Borel-Cantelli
lemma that for a.e. $x\in X,\exists\ k_x\ge 1$ such that
$\b_k(x_k)\le m_k^{{3\over 4}}\ \forall\ k\ge k_x.$
\par It follows from (3) that
$$\sum_{k=1}^\infty m(\{x\in X:x_k\ge a_k-2r_k4^{m_k}\})<\infty,$$
whence  by the Borel-Cantelli lemma, for a.e. $x\in X,\exists\ K_x\ge 1$ such that
$(Sx)_k=x_k+r_kN_k\ \forall\ k\ge K_x.$
\par It follows that a.e. $x\in X$, for $k\ge k_x,K_x$,
$$\align |\b_k((Sx)_k)-c\b_k(x_k)|& =(c-e^{{r_k\over\nu_k}})|\b_k(x_k)|\\ &\le
{cm_k^{{3\over 4}}\over\nu_k} ={c\mu_k^{{3\over 4}}\over\nu_k^{{1\over 4}}}\endalign$$
whence by (5),
$$\sum_{k=1}^\infty |\b_k((Sx)_k)-c\b_k(x_k)|<\infty\text{ for  a.e. }x\in X$$
 and $\v\circ S-c\v$ is a coboundary, being
a product type cocycle, the sum of whose partial transfer functions converges.
\heading \S7 Smooth completely squashable $\R$-extensions\endheading
In this section, we construct smooth completely squashable
$\R$-extensions of rotations of the circle. Ergodicity is established
again using proposition 6.2.

\proclaim{Theorem 7.1} Let $T: \,\,\,x\mapsto x+\alpha$ mod 1 be an
irrational circle rotation (the circle
is represented as the unit interval $[0,1)$).
Let $\alpha$ have unbounded partial quotients. Then there exists a real
smooth ergodic
cocycle $F$ such that for every $c\neq 0$ there exists a rotation $S$
for which
$$
  f\circ S-c\cdot F \tag1
$$
is a ($T$-)coboundary:
\smallskip

1. If $1\leq p<\infty$ is an integer and
$\limsup_{n\to\infty} q_{n+1}/q_n^p =\infty$, the cocycle $F$ can be
found in $\Cal C^p$. (For $p=1$ the condition means unbounded partial
quotients.)
\smallskip

2. If $\limsup_{n\to\infty} q_{n+1}/q_n^p =\infty$ for all positive
integers $p$, $F$ can be found in $\Cal C^\infty$.

\endproclaim
It is known (cf. \cite{Ba-Me}) that if $\limsup_{n\to\infty}
q_{n+1}/q_n^p <\infty$, every $F\in\Cal C_0^p$ is a coboundary.

\demo{Proof} We shall prove Statement 1; then, the other one can be derived
rather easily.

Let $a_n$ be the partial quotients, $q_n$ the convergents.
We define

$c_k=k^3e^{k}$,

$r_k=k^4e^{k}$, $k=1,2,\dots$.\newline
Since $\limsup_{n\to\infty} q_{n+1}/q_n^p =\infty$, there
exist sequences of positive integers

$\ell_k=2[(a_{n_k}-1)/2r_k]$ (where $[x]$ denotes the integer part of $x$),

$\ell'_k=\ell_k/2$,

$d_k=e^kk^6/q_{n_k}$,

$\bar
d_k=d_k(r_kq_{n_k-1})^p=e^{k(p+1)}k^{4p+6}q_{n_k-1}^p/q_{n_k}$,\newline
where
$$
  \sum_{k=1}^\infty \bar d_ke^k <\infty. \tag{2}
$$

Let us suppose that $n=n_k$ is odd (the case with $n_k$ even is
similar).\newline From the continued fraction expansion we get two
Rohlin towers:

$[\{j\alpha\},\{(q_{n-1}+j)\alpha\})$, $j=0,\dots,q_n-1$ and

$[\{q_n\alpha\},1),[\{(j+q_n)\alpha\},|[j\alpha\})$,
$j=1,\dots,q_{n-1}-1$ \newline
(for $0\leq x<1-\|q_{n-1}\alpha\|$ we thus have
$T^{q_{n-1}}x=x+\|q_{n-1}\alpha\|$).
Let us denote

$I_0=[0,\|q_{n-1}\alpha\|)$ $I_i=T^iI_0$, $i=1,\dots,q_n-1$.

For $j=0,\dots,q_{n-1}-1$, the intervals

$I_{j+q_{n-1}},I_{j+2q_{n-1}},\dots I_{j+a_nq_{n-1}}$
\newline
are adjacent.

For $j=0,\dots,r_k-1$ and $u=0,\dots,q_{n_k-1}-1$ we define
$$\align
  J'_{0,0}&=\underset i=0 \to{\overset \ell_k'-1
  \to{\bigcup}} I_{i\cdot q_{n_k-1}},\qquad J'_{0,j}=T^{j\cdot
  \ell_k\cdot q_{n_k-1}}J'_{0,0},\\
  J''_{0,j}&=T^{\ell_k'\cdot
  q_{n_k-1}}J'_{0,j},\\
  J'_{u,j}&=T^uJ'_{0,j},\\
  J''_{u,j}&=T^uJ''_{0,j},\\
  J_{u,j}&=J'_{u,j}\cup J''_{u,j},\\
  J_u&=\underset j=0 \to{\overset r_k-1 \to{\bigcup}}
  J_{u,j}.\endalign
$$
Notice that for every $u$ the sets $J_{u,j}$ are adjacent intervals
composing the interval
$J_u=[\{u\alpha\},r_k\cdot\ell_k\cdot\|q_{n_k-1}\alpha\|)$; each
$J_{u,j}$ is cut in the middle into $J'_{u,j}$ and $J''_{u,j}$. \newline
Let $\bar F_k$ be a $\Cal C^\infty$ function on $[0,1)$ which is

- zero out of $J'_{0,0}$,

- $d_k$ on the middle half of $J'_{0,0}$ \par (i.e. on the
interval
$[\ell_k'\|q_{n_k-1}\alpha\|/4,3\ell_k'\|q_{n_k-1}\alpha\|/4]$),

- has values between 0 and $d_k$ on the rest of $J'_{0,0}$,

$F^{(i)}(0)=0=F^{(i)}(\ell_k'\|q_{n_k-1}\alpha\|)$ for $i=1,\dots,p$.\newline
Moreover, the functions $\bar F_k$ can be found such that there exists
a positive constant $C$,
$$
  \|\bar F_k\|_{\Cal C^p} <C \bar d_k
$$
for all $k$. Let us show this for $p=1$:\newline
Let $\bar f$ be a function which is zero on

$[0,1)\setminus
((0,\ell_k'\|q_{n_k-1}\alpha\|/4)\cup (3\ell_k'\|q_{n_k-1}\alpha\|/4,
\ell_k'\|q_{n_k-1}\alpha\|)),$\newline
on the interval

$(0,\ell_k'\|q_{n_k-1}\alpha\|/4)$\newline
it is a tent-like function, and on the interval

$(3\ell_k'\|q_{n_k-1}\alpha\|/4,
\ell_k'\|q_{n_k-1}\alpha\|)$\newline
it is a reversed tent-like function, both of height
$(8d_k/(\ell_k'\|q_{n_k-1}\alpha\|))$. The indefinite integral $\bar
F_k(t)=\int_0^t \bar f(x)\,dx$ will thus be zero on $[0,1)\setminus
J'_{0,0}$, $d_k$ on the middle half of $J'_{0,0}$ and monotone on the
remaining part of $J'_{0,0}$; there exists a constant $C$ such that
$\|\bar F_k\|_{\Cal C^1}/\bar d_k <C$ for every $k$.\newline
The case of larger exponents $p$ is left to the reader (it can be done
in a recursive way; the constant $C$ depends on $p$).

We
define $\tilde{F}_k=\bar F_k-\bar F_k\circ T^{-\ell_k'\cdot q_{n_k-1}}$
(i.e. $\tilde{F}_k=\bar F_k$ on $J'_{0,0}$; on $J''_{0,0}$, $\tilde{F}_k$
is got by shifting $\bar F_k$ by $\ell_k'\|q_{n-1}\alpha\|$ and changing
the sign);
$$\alignat3
  &F_k=(-1)^j(1+\frac1{c_k})^j\tilde{F}_k\circ T^{-(u+j\cdot \ell_k\cdot
  q_{n_k-1})}
  \,\,\, &&\text{on}\,\,\,
  J_{u,j},\,\,\,&&j=0,\dots,r_k-1,\\
  & && &&u=0,\dots,q_{n_k-1}-1\\
  &F_k=\qquad\qquad 0 &&\text{otherwise}.\endalignat
$$
We have
$$\gather
  |F_k|\leq d_k\cdot (1+\frac1{c_k})^{r_k}\leq d_k\cdot e^k,\\
  \|F_k\|_{\Cal C^p}\leq C\cdot \bar d_k\cdot e^k. \endgather
$$
By (2) we have $\|F_k\|_{\Cal C^p}<\infty$, hence for each subset $K$
of $\Bbb N$ there exists a $\Cal C_0^p$ function
$$
  F_{(K)}=\sum_{k\in K} F_k.
$$

Let  $0\leq j\leq r_k-1$ and $x\in I_{j\cdot\ell_k\cdot
q_{n_k-1}}$.\newline
For $i=u+p\cdot q_{n_k-1}$ where $0\leq u\leq q_{n_k-1}-1$, $0\leq
p\leq\ell'_k-1$ we have $T^ix\in J'_{u,j}$, and for $i=u+p\cdot
q_{n_k-1}+\ell'_k\cdot q_{n_k-1}$ ($0\leq u\leq q_{n_k-1}-1$, $0\leq
p\leq\ell'_k-1$) we have $T^ix\in J''_{u,j}$.\newline
>From the definition of $F_k$ we thus get
$$
  \sum_{i=0}^{\ell'_kq_{n_k-1}-1}
  F_k(T^ix)=-\sum_{i=\ell'_kq_{n_k-1}}^{\ell_kq_{n_k-1}-1} F_k(T^ix)
  \qquad (x\in I_{j\cdot\ell_k\cdot q_{n_k-1}}), \tag3
$$
hence for every $0\leq j\leq r_k-2$ and $x\in I_0$,
$$
  \sum_{i=j\cdot\ell_k\cdot q_{n_k-1}}^{(j+1)\cdot\ell_k\cdot
  q_{n_k-1}-1} F_k(T^ix)=0 \tag{3'}
$$
and
$$
  F_k=G_k-G_k\circ T \tag4
$$
where
$$\alignat2
  G_k(T^ux)&=-\sum_{i=0}^{u-1} F_k(T^ix) \quad &&\text{for}\,\,\,\,x\in
  I_0,\,\,\,\,u=0,\dots,q_{n_k-1}\\
  G_k(x)&=0 &&\text{for}\,\,\,\,x\in [0,1)\setminus \underset i=0
  \to{\overset q_{n_k}-1 \to{\bigcup}} I_i.
  \endalignat
$$
Therefore, $F_k$ is a coboundary with the transfer function $G_k\in \Cal
C^p$.\newline
Let us compute $\sup |G_k|$. We have
$$
  \sup |G_k|=\sup_{x\in I_0} \max\{|\sum_{i=0}^{u-1} F_k(T^ix)|:\,
  0\leq u\leq q_{n_k}-1\};
$$
by (3'), the partial sums are zero for every $u=j\cdot\ell_k\cdot
q_{n_k-1}$, $1\leq j\leq r_k-$. From this and from (3) we get
$$\gathered
  |G_k|\leq \sup \{|\sum_{i=0}^{u-1} F_k\circ
  T^i(x)|:\,x\in X,\ 1\leq u\leq \ell_k\cdot q_{n_k-1}\}\leq \\
  \ell_k\cdot q_{n_k-1} \sup |F_k| =\\
  \ell_k\cdot q_{n_k-1}\cdot  C\cdot
  d_k\cdot(1+\frac1{c_k})^{r_k}
  \leq C\cdot e^k\cdot k^2. \endgathered
  \tag5
$$
As $T^{q_{n_k}-1}$ is the shift (mod 1) by $\|q_{n_k-1}\alpha\|$,
$T^{j\cdot q_{n_k-1}}$,
$j=1,\dots,[\ell'_k/k]$, $k=1,2,\dots$, is a rigid
time. For any fixed positive integer $p$ we thus have
$$
  \lim_{k\to\infty}\max_{j=1,\dots,[\ell'_k/k]}
  |S_{j\cdot q_{n_k-1}}(\sum_{i=0}^p F_i)|=0.
$$
>From this and from $|F_k|\leq d_k\cdot e^k\to 0$ (cf. (2)) follows that
there exists an infinite subset $K\subset \Bbb N$ s.t.
$$
  \lim_{k\in K,k\to\infty}\max_{j=1,\dots,[\ell'_k/k]}
  |S_{j\cdot q_{n_k-1}}(F_{(K)}-F_k)|=0        \tag6
$$
where $F_{(K)}=\sum_{k\in K}F_k$. The set $K$ can be chosen such that
$$
  \sum_{k\in K} \frac1{k} <\infty .\tag7
$$
Let $\Cal A_k$ be the partition of $[0,1)$ into the sets $J_{u,j}$ and
the complement of their union; as $a_{n_k}\to\infty$,
$\Cal A_k\nearrow \Cal A$ (for a subsequence of the
numbers $k$).\newline

Let $0\leq u\leq q_{n_k-1}-1$, $0\leq j\leq r_k-1$ be fixed, $E$ be the
middle third of the interval $J_{u,j}'$. If $1\leq i\leq \ell_k'/k$ and
$k\geq 12$,
$$
  S_{i\cdot q_{n_k-1}}(F_k)=i\cdot q_{n_k-1}\cdot
  (-1)^j(1+\frac1{c_k})^jd_k.
$$
on $E$.\newline
Let $a$ be a fixed number of the same sign as $(-1)^j$ and let
$\epsilon>0$. Without loss of generality we can suppose that
$j$ is even, $a\geq 0$.
By (2) we have
$$
  q_{n_k-1}d_k(1+\frac1{c_k})^{r_k}\leq
  q_{n_k-1}d_k e^{k}\to 0
$$
and for $k$ sufficiently big,
$$
  d_kq_{n_k-1}[\ell'_k/k]\geq
  \frac{e^kk^6}{q_{n_k}}\cdot q_{n_k-1}(\frac{a_{n_k}}{2k^4e^k}-1)\frac1{k}
  \geq k/3,
$$
hence if $k$ is bigger than some constant $k(a,\epsilon)$, then there
exists $1\leq i\leq \ell_k'/k$ for which $S_{i\cdot q_{n_k-1}}(F_k)
\in \Cal U_\epsilon(a)$ on $E$.\newline
The rotation $T^{q_{n_k-1}}$ is the shift by $\|q_{n_k-1}\alpha\|$ (mod
1), hence, if $k\geq 30$, we have
$$
  \lambda(E\cap T^{-i\cdot q_{n_k-1}}E)>0,9\,\lambda(E)
$$
for every $1\leq i\leq \ell_k'/k$.\newline
>From this and from $\lambda(E)\geq \lambda(J_{u,j})/6$ we get the rigid EVC
for $F_k$.
>From the rigid EVC for $F_k$ and (6) follows that for any infinite subset $K'$
of $K$,
$F_{(K')}=\sum_{k\in K'}F_k$ satisfies the essential value condition.
\newline
This way we have found an uncountable set of ergodic cocycles $F=
F_{(K')}\in \Cal C_0^p$.
\medskip

It remains to prove that the set $K$ can be chosen so that for every
$c\neq0$ there exists a rotation $S$ for which
$$
  F\circ S-c\cdot F\qquad\text{is a}\,\,\, T-\text{coboundary}. \tag1
$$
If $F\circ S'-c'\cdot F$ and $F\circ S''-c''\cdot F$ are coboundaries,
then $F\circ S'\circ S''-c'\cdot c''\cdot F =F\circ S'\circ S''
-c''\cdot (F\circ S') +c''\cdot (F\circ S'-c')$ is a coboundary, too,
hence the set of numbers $c$ for which (1) holds true is a group. It
thus suffices to find $S$ for $|c|>1$. First we shall show the proof for
$c>1$. Let us suppose that the $c>1$ is fixed.

Let $j(k)$ be the greatest positive even integer for which
$$
  (1+\frac1{c_k})^{j(k)}<c.
  $$

For nonnegative integers $k$, $v$ let us define a number

$\sigma(k,v)=\{(v(k)+j(k)\cdot\ell_kq_{n_k-1})\alpha\}$ \newline
and a rotation

$\sigma(k,v)=\{(v(k)+j(k)\cdot\ell_kq_{n_k-1})\alpha\}$\newline
(we denote both by the same symbol).\newline
We'll recursively define $K'=\{k_0<k_1<\dots\}\subset K$, nonegative
integers $\{v(k):\,k\in K'\}$, $k=k_0,k_1,\dots$, numbers and rotations
$\{\sigma(k):\,k\in K'\}$ (denoted by the same symbol):\newline
For $k_0$ we choose the smallest element of $K$ and define

$v(k_0)=0$, $\sigma(k_0)=\sigma(k_0,v(k_0))$.\newline
If $k_i$, $v(k_i)$ have been defined for $i=0,\dots,m$, we define
$k_{m+1}$ as the smallest $k\in K$ such that:

1. $k>k_m$.

2. There exists an integer $0\leq v=v(k) <q_{n_k-1}/k$ such that
$$\gather
  |\sigma(k_m)-\sigma(k,v)|<1/2^m,\\
  \sup\,|G_j\circ\sigma(k_m)-G_j\circ\sigma(k,v)|<1/2^m\quad
  \text{for}\,\,\,\,j=k_0,\dots,k_m.
  \endgather
$$
Set $\sigma(k_{m+1})=\sigma(k_{m+1},v(k_{m+1}))$.\newline
The numbers $\sigma(k_m)$ then for $m\to\infty$ converge to a limit
$\sigma$.
\newline
By $S$ we denote the rotation $x\mapsto x+\sigma$ mod 1,
$K'=\{k_0,\,k_1,\dots\}$.\newline
For $k=k_m\in K'$
$$\multline
  \|G_k\circ S-G_k\circ\sigma(k)\|_\infty\leq\\
  \sum_{i=0}^\infty
  \|G_k\circ\sigma(k_{m+i})-G_k\circ\sigma(k_{m+i+1})\|_\infty
  <\sum_{i=0}^\infty 1/2^{m+i} =1/2^{m-1},
  \endmultline
$$
hence
$$
  \sum_{k\in K'} \|G_k\circ S-G_k\circ\sigma(k)\|_\infty
$$
converges.

We have
$$\gather
  F\circ S -c\cdot F=
  \sum_{k\in K'} (F_k\circ S-c\cdot F_k)=\\
  \sum_{k\in K'} \big((F_k\circ S-F_k\circ \sigma(k))+(F_k\circ
  \sigma(k)-(1+\frac1{c_k})^{j(k)}F_k)+((1+\frac1{c_k})^{j(k)}-c)F_k\big).
  \endgather
$$
By (4), each of the functions $F_k$ is a coboundary, hence all
summands in the last sum are coboundaries, too. For proving (1) it
suffices to show that the sum of the corresponding transfer functions
converges:
\smallskip

1. $\sum_{k\in K'} (F_k\circ S-F_k\circ \sigma(k))$.\newline
Every $F_k$ is a coboundary with a transfer
function $G_k$. We have shown that

$\sum_{k\in K'}(G_k\circ S-G_k\circ \sigma(k))$\newline
converges; it is a transfer function of
$\sum_{k\in K'} (F_k\circ S-F_k\circ \sigma(k))$.
\smallskip

2. $\sum_{k\in K'} (F_k\circ
\sigma(k)-(1+\frac1{c_k})^{j(k)}F_k)$.\newline
The function $F_k\circ \sigma(k)-(1+\frac1{c_k})^{j(k)}F_k$ is a
coboundary with a transfer function
$$
  \tilde{G_k}=G_k\circ \sigma(k)-(1+\frac1{c_k})^{j(k)}G_k.
$$
>From the definitions of $F_k$, $G_k$, and $\sigma(k)$ follows that for
$x\in I_i$,
$i=u+t\cdot \ell_k\cdot q_{n_k-1}+p\cdot q_{n_k-1}$, $0\leq u$,
$u+v(k)\leq q_{n_k-1}$, $0\leq t$, $t+j(k)\leq r_k-1$, and $0\leq p\leq
\ell_k-1$, we have
$$
  \tilde{G_k}=G_k\circ \sigma(k)-(1+\frac1{c_k})^{j(k)}G_k=0.
$$
>From $v(k)\leq q_{n_k-1}/k$ and $j(k)\leq c_k\cdot(\log c+1)=(\log
c+1)\cdot r_k/k$ (for $k$ sufficiently big) we get that
$$
  \lambda(\tilde{G}_k\neq 0)<(3+\log c)/k
$$
for $k$ big enough. From this and (7) follows that  the sum
$\sum_{k\in K'} \tilde{G}_k$  converges almost surely.
\smallskip

3. $\sum_{k\in K'} ((1+\frac1{c_k})^{j(k)}-c)F_k$.\newline
By (4), $F_k$ is a
coboundary with a transfer function $G_k$ and by (5), $|G_k|$ is bounded by
$(1/2)e^kk^2$.
We can easily see that $c-(1+\frac1{c_k})^{j(k)}\leq 2c/c_k$, hence
$$
  \sum_{k\in K'} |((1+\frac1{c_k})^{j(k)}-c)G_k|\leq \sum_{k\in K'}
  c\cdot e^kk^2/c_k\leq \sum_{k\in K'} 1/k
$$
where the last sum is  finite by (7).
\smallskip

If we define $j(k)$ as the biggest odd number for which
$$
  (1+\frac1{c_k})^{j(k)}<c,
$$
we get the rotation $S$ for $c<-1$, hence (1) holds true for all $c\neq
0$.
\smallskip

All the cocycles $F_k$ which we defined are from $\Cal C_0^\infty$. If
$\limsup_{n_k\to\infty} q_{n_k+1}/q_n^p =\infty$, we can find the cocycles
$F_k$ such that $\sum_{k=1}^\infty F_k$ converges in every $\Cal C^p$,
$1\leq p<\infty$, hence $F=\sum_{k=1}^\infty F_k\in\Cal C_0^\infty$. This
proves the second statement of the theorem.
\hfill\qed
\enddemo

\

\heading References\endheading

\enddocument